\documentclass[pdflatex, sn-basic]{sn-jnl}       

\usepackage{tikz}
\usetikzlibrary {math}

\usepackage{graphicx}%
\usepackage{multirow}%
\usepackage{amsmath,amssymb,amsfonts}%
\usepackage{amsthm}%
\usepackage{mathrsfs}%
\usepackage[title]{appendix}%
\usepackage{xcolor}%
\usepackage{textcomp}%
\usepackage{manyfoot}%
\usepackage{booktabs}%
\usepackage{algorithm}%
\usepackage{algorithmicx}%
\usepackage{algpseudocode}%
\usepackage{listings}%
\usepackage{mwe}%
\usepackage{float}%
\usepackage{dirtytalk}%
\usepackage{svg}%
\usepackage{makecell}%

\renewcommand{\orcidlogo}{%
  \includegraphics[width=10pt]{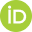}%
}
\renewcommand{\orcid}[1]{\href{https://orcid.org/#1}{\orcidlogo}}

\usepackage{todonotes} 
\usepackage{anyfontsize} 

\graphicspath{{figures/}}

\RequirePackage{fix-cm}
\begin{document}
 \newcommand{\be}{\begin{equation}}
\newcommand{\ee}{\end{equation}}
\newcommand\cC{{\cal C} }
\newcommand\cF{{\cal F} }
\newcommand\cG{{\cal G} }
\newcommand\Fb{{\overline{\! F} }}
\newcommand\Hb{{\overline{\! H} }}
\newcommand\rH{{\rm H }}
\newcommand\rHh{\widehat{\rm H }}
\newcommand\cH{{\cal H} }
\newcommand\bI{{\bf I} }
\def\N{\mathbb{N}}
\newcommand\gO{{\cal O} }
\newcommand\pt{\tilde{p} }

\def\R{\mathbb{R}}
\newcommand\br{{\bf r} }
\newcommand\brt{\tilde{\bf r} }
\newcommand\rt{\tilde{ r} }
\newcommand\xb{{\overline x} }
\newcommand\xt{{\tilde x}}
\newcommand\bx  {{\boldsymbol{x} }}
\newcommand\bxt  {\tilde{\boldsymbol{x} }}
\def\Z{\mathbb{Z}}
\newcommand\bZ {{\boldsymbol{Z} }}

\newcommand\lam{{\lambda} }
\newcommand\lat{{\tilde\lambda} }
\newcommand\Lam{{\Lambda} }
\newcommand\gam{{\gamma} }
\newcommand\eps{{\varepsilon} }
\newcommand\Om{{\Omega} }
\newcommand\om{{\omega} }
\newcommand\zetat{{\tilde \zeta} }
\newcommand\bzeta{\boldsymbol{\zeta}}
\newcommand\nut{\tilde\nu}
\newcommand\bzero{\boldsymbol{0}}

\newcommand{\Las}{\Lam^{\star} }
\newcommand{\as}{a^{\star} }
\newcommand{\ns}{n^{\star} }

 \newcommand{\norm}[1]{\left\Vert#1\right\Vert}
 \newcommand{\qtext}[1]{\quad \text{#1}\quad}

\newcommand{\der}[2]{\frac{d#1 \hfill}{d#2 \hfill}}
\newcommand{\dron}[2]{\frac{\partial#1 \hfill}{\partial#2 \hfill}}

\newcommand{\FL}[1]{\cF_{L_#1} }
\newcommand{\VFL}{${\cal VF}_L$ }
\newcommand{\VFE}{${\cal VF}_E$ }

\newcommand{\mnras}{MNRAS }

\newcommand{\tpr}[1]{\textcolor{red}{#1}}
\newcommand{\tap}[1]{\textcolor{blue}{#1}}
\newcommand{\tme}[1]{\textcolor{green}{#1}}


\newcommand\pdouze{$P_{12}$}
\newcommand{\rad}{{\rm rad}}
\newcommand{\VFLtitre}{\texorpdfstring{\VFL}{VFL}}
\newcommand{\changed}[1]{#1}
\title[Jacobi]{Marchal's family of periodic orbits. I: Stability of inclined co-orbital planetary systems}
\author*[1]{\fnm{Alexandre} \sur{Prieur} \orcid{0009-0004-9125-4678}}\email{alexandre.prieur@obspm.fr}

\author[1]{\fnm{Philippe} \sur{Robutel} \orcid{0000-0001-8932-1713}}\email{philippe.robutel@obspm.fr}


\affil[1]{\orgdiv{LTE}, \orgname{Observatoire de Paris -- PSL}, \orgname{CNRS}, \orgname{Sorbonne Université}}

\date{\today}


\setcounter{tocdepth}{3}

\abstract{At the Lagrange relative equilibrium of the three-body problem, for all values of the masses, the elliptic eigenvalues associated with vertical eigenvectors give rise to spatial quasi-periodic orbits, which become periodic in a rotating frame.
In 2009, by averaging out the fast frequencies, Christian Marchal showed that these orbits, which are fixed points in the restricted average problem, form a one-parameter family connecting $L_4$ to $L_5$.
Using perturbation methods, we show the persistence of this family in the average three-body problem for nonzero masses in the limit where one mass is dominant over the other two (known as the planetary problem).
We also give an analytical approximation valid for mutual inclinations less than $60^\circ$.
Then, using purely numerical methods, we show that this family exists in the full three-body problem (neither restricted nor average) for a wide range of masses, beyond the planetary case. We also show that the stability of its orbits evolves along the family, with inclined systems remaining stable for masses exceeding the Gascheau’s value (also known as Routh’s critical value).
Finally, we show the impact of this family's stability on the global dynamics of the co-orbital region as well as its high instability for mutual inclinations exceeding $60^\circ$.
}
\keywords{Three-body problem, Planetary problem, Numerical methods, Periodic orbits, Orbit stability}



\maketitle

\section{Introduction}
\label{sec:intro}

The only known explicit solutions of the 3-body problem correspond to two specific configurations: the aligned configurations discovered by Euler (\citeyear{Euler1764}) and the equilateral configurations demonstrated by Lagrange (\citeyear{Lagrange1772}), where each of the three bodies follows a Keplerian motion.
The question of the spectral stability of these configurations was only resolved in the middle of the XIX\textsuperscript{th} century. Liouville showed in (\citeyear{Liouville1842}) that, contrary to what Euler seemed to think, the aligned configurations are unstable regardless of the masses of the bodies. The following year, \cite{Ga1843} addressed the problem of the spectral stability of equilateral configurations and demonstrated that the motion is either stable or unstable, depending on whether the ratio of two terms is greater or less than 27: the square of the sum of the three masses, and the sum of the pairwise products of these masses, \changed{or $(\sum_i m_i)^2 / \sum_{i < j} m_i m_j$}. Thus, spectral stability occurs only when one of the masses is predominant compared to the other two.
Gascheau established this criterion only in the case where the three bodies describe a circular motion (relative equilibrium). For the elliptical case (central configuration), one may refer to  \cite{Dan1964} or \cite{Bennett1965} in the restricted case (one of the masses equal to zero) and  to \cite{Nauenberg2002, Robe2002, MR3218836} in the general case.

To illustrate Gascheau's criterion with a situation that we will focus on in this paper, let us assume that two of the masses are equal to one another, say equal to $\eps$, and that the sum of the three masses is equal to 1. In this case, the equilateral equilibrium, in the case of circular motion, is linearly stable as long as $\eps$  is lower than $(3 - 2\sqrt2)/9 \approx 0.0191$. In the following, we will refer to this quantity as the \say{Gascheau value}. Let us note that, in the case of the restricted problem, this critical value is equal to $(1 - \sqrt{69/18})/2 \approx 0.0385$.

Below the Gascheau value, the periodic orbit corresponding to the circular equilateral equilibrium -- or the fixed point if placed in a reference frame rotating at the correct frequency -- is elliptic. Thus, from this equilibrium, a Lyapunov family of periodic orbits arises for each pair of conjugate eigenvalues \cite[see][]{MeHa1992}.
Three families corresponding to coplanar orbits emanate from the equilateral relative equilibrium. Among these families, the most known is that of so-called homographic (in the sense of self-similar) orbits, where each body follows a Keplerian orbit with the same eccentricity. It starts with zero eccentricity and ends with an eccentricity of 1, leading to a triple collision. A comprehensive study of these orbits and their stability can be found in \cite{Robe2002}.
A second family corresponds to nearly circular orbits with eccentricities of order $\sqrt{\eps}$, the parameter $\eps$ measuring the smallness of the masses. It is worth mentioning that, in the restricted problem, this family coincides with the so-called \say{long period family}. The last family, called anti-Lagrange by \cite{GiuBeMiFe2010}, which does not exist in the restricted problem, consists of eccentric orbits whose pericenters precess at a common frequency of order $\eps$. A study of these three families is made in \cite{RoPo2013}.

\changed{Although spatial co-orbital motion has been studied in other cases (see e.g. \cite{voyatzisQuasisatellitePeriodicMotion2018} for a family of quasi-periodic orbits emerging from planar quasi-satellite)}, the periodic orbits (in rotating frame) that emerge from the Lagrange relative equilibrium in the vertical direction have been much less investigated. As far as we know, precise results have been obtained only in two extreme cases: the three equal masses, and the restricted problem involving the Sun, Jupiter, and a test particle.

In the first case, \cite{Marchal1990, Marchal2000} conjectured the existence of the \pdouze\, family, a family of periodic orbits in a rotating frame originating from the equilateral relative equilibrium  in the vertical direction (existence of mutual inclination between the orbits) to join the choreographic figure-Eight orbit of \cite{CheMo2000}. All orbits have the same symmetry group of order 12 as the figure-eight. The local existence and uniqueness of the \pdouze\,  family near the Lagrange relative equilibrium was proven by \cite{FeChe2008}, while the existence of the whole family was demonstrated by \cite{Calleja_etal2024}, with computer assistance.

In contrast, in the case of the restricted problem, \cite{Marchal2009} demonstrates, through averaging over the fast component of the motion, that there exists a one-parameter family of fixed points of the average problem (this parameter can be the inclination of the Trojan's orbit relative to that of Jupiter) connecting $L_4$ to $L_5$.
Along this family, the difference between the mean longitude of the asteroid and that of Jupiter (which is constant for a given orbit) evolves continuously between  $60^\circ$ and $300^\circ$. Likewise, along this path, the inclination of the test particle evolves: it is equal to zero at $L_4$ and $L_5$ and reaches a maximum close to $146^\circ$ degrees when the difference in longitude is $180^\circ$.
If we return to an inertial reference frame, these orbits are characterized by two frequencies: a slow frequency (of the same order of magnitude as the ratio of Jupiter's mass to that of the Sun) associated with the precession of the ascending node, and a fast frequency over which the averaging has been performed. They can also be viewed as periodic orbits in a frame rotating with the node.

At the end of his paper, Marchal conjectures that the orbits described above belong to a three-parameter family of periodic orbits (in a rotating reference frame):

\begin{quote}
    These periodic orbits belong to a three parameter family of similar periodic orbits of the general three body problem, two parameters being the mass ratios of the three bodies and the third a generalization of the inclination. This family of periodic orbits always begins at the circular Lagrangian motion of the three masses of interest and “escapes” in the normal direction. The one parameter family P12 of \cite{Marchal2000} is another limit case of this three parameter family: the case of three equal masses.
\end{quote}

In the following sections, we will focus on a part of the conjecture stated above, where the associated families of orbits will be referred to as {\it Marchal family}.
In particular, we will attempt to address the following questions: Does the family found by Marchal exist in the complete problem? Does it persist for non-zero masses? What is the stability of the orbits along this family? Does this have an influence on the global stability of the co-orbital resonance?

Although \cite{Leleu2017PhD} found numerical hints indicating that the family probably exists in the full problem, this topic does not yet seem to have been addressed precisely.
In order to address these questions, we will consider the planetary problem where two of the bodies have a small non-zero mass compared to the third. We will first follow the steps of \cite{Marchal2009} by focusing, in section \ref{sec:moyen}, on the average 3-body problem, which is relevant due to the perturbative regime.
After reducing the problem (section \ref{sec:Jacobi}), we will show the existence of a family of fixed points forming a one-parameter family similar to that described by Marchal in the case of the restricted problem. This will provide a generalization of Marchal's result in the case of the average planetary problem. Then, we will address the complete (i.e. non-average) problem in section \ref{sec:num}; after introducing a new set of coordinates, we will present numerical evidence of Marchal family's existence in the complete, non-restricted problem. We will finally show the family's stability limits, from the planetary problem to beyond Gascheau's value.

\section{The average co-orbital problem} 
\label{sec:moyen}
\subsection{Reduction and averaging of the planetary Hamiltonian}

\subsubsection{Jacobi's reduction}
\label{sec:Jacobi}

\begin{figure}[!ht]
%
%

\tikzset{math3d/.style={x= {(-0.353cm,-0.353cm)}, z={(0cm,1cm)},y={(1cm,0cm)}}}
\tikzstyle{orbiteA} = [color=red,font=\scriptsize]
\tikzstyle{orbiteB} = [color=blue,font=\scriptsize]
\tikzstyle{ref} = [color=black,font={\scriptsize}]
\tikzstyle{Noeud} = [color=black,font=\scriptsize]

\tikzmath{
\Inc  = 30;
\IncB  = 30;
\cI = cos(\Inc);
\sI = sin(\Inc);
\r = 3;
\rfA = 3.1;
\rfB = 2.7;
\RotAxe= -22;
\cRotAxe = cos(\RotAxe);
\sRotAxe = sin(\RotAxe);
}

\begin{tikzpicture}[math3d,scale=1.45]

\draw [domain=pi:2*pi, samples=80, smooth,orbiteA] plot ({\r*cos(\x r)}, {\r*\cI*sin(\x r)}, {\r*\sI*sin(\x r)}) ;
\draw [domain=0:pi, samples=80, smooth,orbiteB] plot ({\r*cos(\x r)}, {\r*\cI*sin(\x r)}, {-\r*\sI*sin(\x r)}) ;
%

%
\fill[gray!20,opacity=0.8] (-\r,-\r,0) -- (\r,-\r,0) -- (\r,-\r,0) -- (\r,\r) -- (-\r,\r) -- cycle;
\draw [domain=pi:3*pi/2, samples=80, smooth,orbiteB,dashed,->,>=latex] plot ({\rfA*cos(\x r)}, {\rfA*\cI*sin(\x r)}, {-\rfA*\sI*sin(\x r)}) ;
\draw [domain=-.05:pi/2*1.17, samples=80, smooth,orbiteA,dashed,->,>=latex] plot ({\rfB*cos(\x r)}, {\rfB*\cI*sin(\x r)}, {\rfB*\sI*sin(\x r)}) ;
\draw [domain=0:pi, samples=80, smooth,orbiteA] plot ({\r*cos(\x r)}, {\r*\cI*sin(\x r)}, {\r*\sI*sin(\x r)}) ;
\draw [domain=pi:2*pi, samples=80, smooth,orbiteB] plot ({\r*cos(\x r)}, {\r*\cI*sin(\x r)}, {-\r*\sI*sin(\x r)}) ;
%
\draw[thick, >=latex,->] (0,0,0) -- (0,0,\r);
\draw[thick,dashed] (-\r*1.4,0,0) -- (\r*1.4,0,0);
\draw[thin,>=latex,->] (0,0,0) -- (\r*\cRotAxe*1.4,\r*\sRotAxe*1.4,0);
\node[ref] at (\r*\cRotAxe*1.43,\r*\sRotAxe*1.43,0) {$x$};
%
\draw[thin,>=latex,->,orbiteA,dashed] (\r*\cRotAxe*.5,\r*\sRotAxe*.5,0) arc (\RotAxe:0:\r*.5);
\draw[thin,>=latex,->,orbiteB,dashed] (\r*\cRotAxe*.6,\r*\sRotAxe*.6,0) arc (\RotAxe:180:\r*.6);
\draw [domain=0:\Inc, samples=80, smooth,orbiteA,->,>=latex,thin,densely dashed] plot ({\r}, {cos(\x)}, {sin(\x)}) ;
\def\coefa{.78};\def\coefb{1.6};
\draw [domain=0:\Inc, samples=80, smooth,orbiteB,->,>=latex,thin,densely dashed]
plot ({-\r}, {-\coefb*cos(\x*\coefa)}, {\coefb*sin(\x*\coefa)}) ;
%
\node[orbiteB,left] at (0,.6*\r) {$\Omega_2$};
\node[orbiteB] at (-\r*1.2,-\r/1.8) {$I_2$};
\node[orbiteA] at (\r*.8,\r*0.3) {$I_1$};
\node[orbiteA] at (\r*.38,-\r*0.1) {$\Omega_1$};
\node[ref] at (-\r*3,-\r*1.05) {${\cal C}$};

%

\def\AnglePa{90} 
\coordinate(Pa) at  ({\r*cos(\AnglePa)}, {\r*\cI*sin(\AnglePa)}, {\r*\sI*sin(\AnglePa)}) ;
\node[orbiteA] at (Pa)  {$\bullet$};
\node[orbiteA,right] at (Pa)  {$P_1$};
\def\AnglePb{-90}
\coordinate(Pb) at ({\r*cos(\AnglePb)}, {\r*\cI*sin(\AnglePb)}, {-\r*\sI*sin(\AnglePb)});
\node[orbiteB] at (Pb)  {$\bullet$};
\node[orbiteB,right] at (Pb)  {$P_2$};
\node[orbiteA,left=0.2cm] at (Pa)  {$w_1$};
\node[orbiteB,above=.4cm] at (Pb)  {$w_2$};

\def\AnglePa{0} 
\coordinate(Pa) at  ({\r*cos(\AnglePa)}, {\r*\cI*sin(\AnglePa)}, {\r*\sI*sin(\AnglePa)}) ;
\node[Noeud] at (Pa)  {$\bullet$};
\node[orbiteA,right] at (Pa)  {$N_1$};
\def\AnglePb{180}
\coordinate(Pb) at ({\r*cos(\AnglePb)}, {\r*\cI*sin(\AnglePb)}, {-\r*\sI*sin(\AnglePb)});
\node[Noeud] at (Pb)  {$\bullet$};
\node[orbiteB,right] at (Pb)  {$N_2$};

\end{tikzpicture}
 \caption{Coordinate system linked to the invariant plane (perpendicular to the total angular momentum and containing the most massive body). The angular positions of the bodies are indicated by the angles $w_j$ (draconic true anomalies) originating at the ascending node ($N_j$) of the corresponding orbit.}
 \label{fig:Jacobi}
\end{figure}
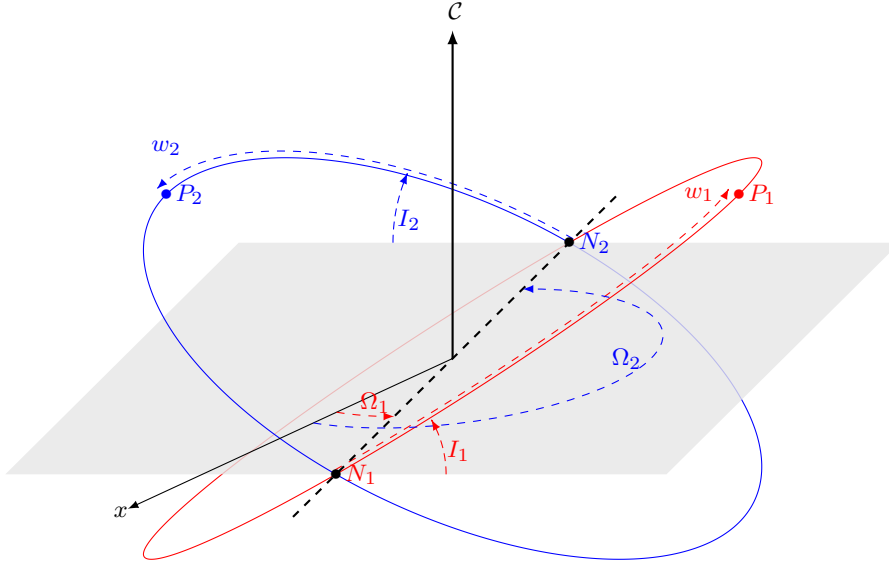

 We consider two planets of respective masses $\eps m_1$ and $\eps m_2$, $\eps$ being a parameter reflecting the smallness of the planetary masses compared to that of the central body of mass $m_0$.
 In heliocentric canonical coordinates \citep{LaRo1995} and rescaling both action variables and time by a factor $\eps$ \citep[see][]{RoNiPo2016}, the Hamiltonian of the planetary 3-body problem reads:

\be
\begin{aligned}
H(\brt,\br) &= H_K(\brt,\br) + \eps H_P(\brt,\br) \quad \text{with}  \\
 H_K(\brt,\br) &=  \sum _{j\in\, \{1,2\}} \left(
 \frac{\brt_j^2}{2\beta_j} - \frac{\mu_j \beta_j}{\norm{\br_j}} 
    \right)
    \quad \text{and} \\
 H_P(\brt,\br) &= \frac{\brt_1\cdot\brt_2}{m_0} - \cG\frac{m_1m_2}{\norm{\br_1 -\br_2}} .
     \end{aligned}
     \label{eq:ham_cart}
\ee
In these expressions,  $\br = (\br_1,\br_2) \in \R^3\times\R^3$ contains the heliocentric positions of planets $1$ and $2$ while $\brt = (\brt_1,\brt_2)$, the conjugated variables of $\br_j$, are associated with the rescaled barycentric linear momentum of the same body. The parameters $\beta_j$ and $\mu_j$ are defined as:
\be
\beta_j = \frac{m_jm_0}{m_0 + \eps m_j} \qtext{and} \mu_j = \cG(m_0 + \eps m_j).
\label{eq:beta}
\ee
The invariance of the problem under rotation around the axis of the total angular momentum allows us to reduce the number of degrees of freedom by eliminating inclinations and ascending nodes -- in other words, we perform the Jacobi reduction \cite[see][]{Ro1995,MaEf2023}. To do this, let us introduce Poincar\'e's variables related to the invariant plane -- the plane perpendicular to the total angular momentum and containing the most massive body (see Fig.\ref{fig:Jacobi}) -- defined as:
\be
\begin{array}{ll}
\Lam_j  =  \beta_j\sqrt{\mu_j a_j}\, ,   & \lat_j = M_j + \omega_j\\
x_j = \sqrt{\Lam_j - G_j}\exp(i\om_j) \, , & \xt_j = -i\xb\\
\Psi_1 = G_1 \cos I_1 + G_2 \cos I_2\, , \quad& \psi_1 = (\Om_1 + \Om_2)/2 \\
\Psi_2 = G_1 \cos I_1 - G_2\cos I_2\, , & \psi_2 = (\Om_1 - \Om_2)/2, \\
\qtext{with} G_j = \Lam_j \sqrt{1-e_j^2} ,& 
\end{array}
\label{eq:var_can}
\ee
where, for $j \in \{1, 2\}$ the planet number, $a_j$ is its semi-major axis, $e_j$ its eccentricity, $I_j$ its inclination, $M_j$ its mean anomaly, $\om_j$ the argument of its pericenter, and $\Om_j$ the longitude of its ascending node.  $\xb$ represents the complex conjugated of $x$ and $i=\sqrt{-1}$.

As in this reference frame the ascending nodes of the two planets are in opposition, $\psi_2 = \pi/2$ and its conjugated variable $\Psi_2$ are ignorable. Similarly, the action $\Psi_1$ being equal to the modulus of the total angular momentum, $\psi_1$ is also ignorable. As a result, the Hamiltonian depends only on $8$ of the variables expressed in the coordinates (\ref{eq:var_can}):  $(\Lam_j,\lat_j,x_j,\xt_j)_{j=1,2}$ -- plus an additional parameter $\Psi_1 = \cC$. Since these latter variables do not involve the angles $\Omega_j$ but only $M_j$ and $\omega_j$, using these coordinates amounts to considering the trajectories in a reference frame rotating with the common line of the nodes of the planets. This rotating frame will henceforth be referred to as the draconic frame.

The ascending nodes being in opposition, the mutual inclination $J$ between the two planets is equal to the sum of the individual inclinations:  $J = I_1 + I_2$. Thus, the modulus of the total angular momentum and the mutual inclination are related by the following expression:
\be
\cC^2 =  G_1^2 + G_2^2 + 2 G_1 G_2 \cos J\, .
\label{eq:C_J}
\ee
Let us finally note that, although the inclinations and longitudes of the nodes have been eliminated by Jacobi's reduction, it is still possible to retrieve them using:
\be
\cC\cos I_1 = G_1 +G_2\cos J, \quad \cC\cos I_2 = G_2 +G_1\cos J.
\label{eq:inclin}
\ee
Similarly, the time derivative of the ascending nodes' longitudes are given by:
\be
\dot\Om_1 = \dot\Om_2 = \dron{\cH}{\cC},
\label{eq:noeud}
\ee
where $\cH$ is the Hamiltonian (\ref{eq:ham_cart}) expressed in the variables given by  (\ref{eq:var_can}).

\subsubsection{The co-orbital resonance}
\label{sec:averaged_ham}

\paragraph{Resonant coordinate system}
\label{sec:coord_res}

In this section, we briefly present an approximation of the planetary problem valid in a neighborhood of the co-orbital resonance (1:1 orbital resonance). For more information, we refer the reader to \cite{RoPo2013} and \cite{NiPoRo2020}.

Given an arbitrary orbital frequency $\ns$, we expand the quantities depending on the semi-major axes of the planets near a common semi-major axis $\as$ related to the mean motion $\ns$ by the relation:
\be
\as = \mu_0^{1/3} {\ns}^{-2/3} \qtext{with} \mu_0 = \cG m_0. 
\ee
It follows that the Keplerian orbital frequencies of the two planets evaluated at $\as$ are $\eps$-close, namely:
\be
n_j = \sqrt{\frac{\mu_j}{{\as}^3}} = \ns\sqrt{\frac{\mu_j}{\mu_0}} =  \ns \left( 1 + \gO(\eps) \right) .
\ee
Since the difference in mean longitudes of the two planets is the natural parameter characterizing co-orbital motion, it is natural to  introduce the canonical coordinate system:
 \be
 (Z_j, \zeta_j, x_j, \xt_j)_{j=1,2}
 \label{eq:coord_res}
 \ee
given by the following expression with $\Las_j = \beta_j\sqrt{\mu_j \as} = \beta_j\mu_j^{1/2}\mu_0^{1/6}{\ns}^{-1/3}$ :
\be
\begin{array}{llll}
Z_1 = \Lam_1 - \Las_1 & \quad & \zeta_1  &=  \lat_1  - \lat_2 + \pi\, ,\\
Z_2 = \Lam_1 + \Lam_2 -\left( \Las_1 + \Las_2 \right), & \quad & \zeta_2  &= \lat_2\, ,                         \\
\end{array}
\ee
the $x_j$ and $\xt_j$ remaining unchanged\footnote{$\pi$ is added to $\lat_1  - \lat_2$ so that we find the  equilateral equilibria  at $\zeta_1=\pm\pi/3 $ and not at $\pi \pm\pi/3 $. Note that the difference of the mean longitudes of the two planets is equal to $\lat_1  - \lat_2 + \pi$.}.
The pairs of variables $(Z_1,Z_2)$, $(\zeta_1,\zeta_2)$, $(x_1,x_2)$ and $(\xt_1,\xt_2)$ will be respectively denoted  $\bZ$, $\bzeta$, $\bx$ and $\bxt$.  
In these coordinates, the Hamiltonian (\ref{eq:ham_cart}) reads:
 \be
 \begin{aligned}
 & \rH(\bZ, \bzeta, \bx, \bxt, \cC) = \rH_K(\bZ) + \eps \rH_P(\bZ, \bzeta, \bx, \bxt, \cC) \\
 & \qtext{with} 
  \rH_K(\bZ) = -\frac{\mu^2_1\beta_1^3}{2(\Las_1 + Z_1)^2} -\frac{\mu^2_2\beta_2^3}{2(\Las_2 + Z_2 - Z_1)^2} \\
  & \qtext{and}  \rH_P(\bZ, \bzeta, \bx, \bxt, \cC) = H_P(\brt,\br).
 \end{aligned}
 \label{eq:H_zeta}
 \ee
 
\paragraph{Truncation in \texorpdfstring{$Z_j$}{Z\_j}}

To  carry out our study, it is sufficient to limit ourselves to a neighborhood of $\Las_j$ of size $\sqrt\eps$ , that is, for: 
\be
 \vert Z_j \vert < c\, \sqrt{\eps},
 \label{borneGamma}
\ee
for some $c>0.$
To this end, we first expand the Keplerian part of the Hamiltonian  at $\bZ = \bzero$ to get:
\be
\begin{aligned}
 \rH_K(\bZ)  =  \,
  &\ns Z_2 - \frac{3{\ns}^{4/3}(m_1+m_2)}{2\mu_0^{2/3}m_1m_2}
                     \left[
                     \left(Z_1 - \frac{m_1}{m_1+m_2}Z_2\right)^2
                     + \frac{m_1m_2}{(m_1+m_2)^2}Z_2^2
                     \right]  \\
                  & + \gO(\eps^{3/2})
 \end{aligned} 
 \label{eq:H_Kep_Z}
 \ee
 while a truncation of the perturbation to the zero degree gives:
\be
\eps\rH_P(\bZ,\bzeta,\bx,\bxt,\cC) = \eps\rH_P(\bzero,\bzeta,\bx,\bxt,\cC) + \gO(\eps^{3/2}).
\ee

\subsubsection{Expansion of the perturbation}
 \label{sec:expansion}
According to formula (\ref{eq:ham_cart}), the perturbation $\rH_P$  is composed of two term: a gravitational part involving the mutual distances $\norm{\br_1 - \br_2}$ and a kinetic component proportional to $\brt_1\cdot\brt_2$.
Following \cite{Ro1995}, the mutual distance can be expressed as:
\be
\Delta^2 = \norm{\br_1 - \br_2}^2 = r_1^2 + r_2^2 -2\br_1\cdot\br_2,
\ee
where, $r_j = \norm{\br_j}$ and
\be 
\frac{\br_1\cdot\br_2}{r_1r_2} =   -\cos w^- + (\cos w^- - \cos w^+) s^2
\label{eq:CosS}
\ee
with
\be
 s = \sin(J/2),\, w^+ = w_1 + w_2, \quad w^- = w_1 - w_2, \quad w_j = v_j + \omega_j,
 \ee
 $v_j$ being the true anomaly of the $j$-th planet and $\omega_j$ its pericenter arguments (see Fig. \ref{fig:Jacobi}).
The kinetic component reads:
\be 
\begin{aligned}
&\frac{\brt_1\cdot\brt_2}{m_0} = \frac{\beta_1\beta_2}{m_0}\sqrt{\frac{\mu_1\mu_2}{a_1a_2}} F_I 
\qtext{with} \\
&F_I  =   -(1-e_1^2)^{-1/2}(1-e_2^2)^{-1/2}\left[S_1S_2 + C_1C_2  - 2s^2C_1C_2 \right],  \\
 &\text{and}\quad          S_j = \sin w_j + e_j\sin \omega_j, \quad
             C_j = \cos w_j + e_j \cos\omega_j.
\end{aligned}
\label{eq:indirect}
\ee
%
Since we only need to calculate $\rH_P(\bzero,\bzeta,\bx,\bxt,\cC)$, the quantities $(a_j,\beta_j,\mu_j)$ will be replaced in  (\ref{eq:CosS}) and (\ref{eq:indirect}) by $(\as,m_j,\mu_0)$. We thus get:
\be
\rH_P(\bzero, \, \bx, \bxt, \cC) = \mu_0^{2/3}{\ns}^{2/3}\frac{m_1m_2}{m_0} 
\left(
F_I  - \frac{\as}{\!\!\Delta^\star}
\right) 
\ee
where 
\be
\Delta^\star = 
\as\left(
\sum_{j\in\{1,2\}} \!\! \left[ \frac{1 - e_j^2}{1 + e_j\cos(w_j -  \omega_j)} \right]^2 \! \!\! +2\cos w^- - 2(\cos w^- - \cos w^+) s^2
\right)^{1/2}\!\!\!\!\!\!.
\ee
Since the quantities $e_j, v_j, \omega_j,\cos w^-$ and $\cos w^+$ are independent of the mutual inclination, they  can be expressed in term of $(\bzeta,\bx,\bxt)$, for instance by usual expansion in power of the eccentricities.
 
The final step is to get rid of the mutual inclination $J$ involved in (\ref{eq:CosS}) and (\ref{eq:indirect})  by means of $s$.
To do this, it is convenient to introduce the parameter $J_0$, or equivalently $s_0 = \sin(J_0/2)$, where $J_0$ is the greatest value of the mutual inclination for a given value of $\cC, \Lam_1$ and $\Lam_2$. Note that it is equal to the value of the mutual inclination when $e_1= e_2 = 0$ for a given value of $\cC$.

 Starting from the expression of the total angular momentum in terms of elliptic elements, that is from  (\ref{eq:var_can}):
\be
\cC^2 =  G_1^2 + G_2^2 + 2G_1G_2\cos J = \Lam_1^2 + \Lam_2^2 + 2\Lam_2\Lam_2\cos J_0,
\ee
and inserting this expression of $s_0$ in the previous formula leads to:
\be
\begin{aligned}
s^2 &= \frac12\left( 1 - \cos J\right) \\
    & =  \frac{4s_0^2 + (G_1 + G_2)^2 -(\Lam_1 + \Lam_2)^2}{4G_1G_2} \\
    & =   \frac{4\Lam_1\Lam_2s_0^2 - 2(\Lam_1+\Lam_2)(\vert x_1\vert^2 + \vert x_2\vert^2) - (\vert x_1\vert^2 + \vert x_2\vert^2) ^2 }{4(\Lam_1 - \vert x_1\vert^2) (\Lam_2 - \vert x_2\vert^2)} \\
    & =   \frac{4\Las_1\Las_2s_0^2 - 2(\Las_1+\Las_2)(\vert x_1\vert^2 + \vert x_2\vert^2) - (\vert x_1\vert^2 + \vert x_2\vert^2) ^2 }{4(\Las_1 - \vert x_1\vert^2) (\Las_2 - \vert x_2\vert^2)} + \gO(\sqrt\eps) .
\end{aligned}
\label{eq:J_s}
\ee
This last expression allows us to get rid of the mutual inclination, which completes the reduction of angular momentum.
We still need to calculate  the explicit expression of the perturbation in the variables $(\bzeta, \bx, \bxt)$.

While $F_I$ can be  easily expanded  in power series of $(\bx, \bxt)$, the expansion of the mutual distance requires more effort.
By adapting \cite{Ro1995} to our case, we split $\Delta^2$ as:
\be
\begin{aligned}
&\Delta^2 =  A + V \qtext{with} V = \Delta^2 - A = \gO(e_j)  \qtext{and}\\
&A = 2a_*^2\left( 1 - \cos\zeta_1 + s_0^2( \cos(\zeta_1) - \cos(\zeta_1 + 2\zeta_2))  \right)
\end{aligned}
\ee
Expanded $V$ in power of $e_j$, or equivalently $(\bx,\bxt)$, using standard expressions from Kepler's problem, and substituting the results into the Taylor expansion
\be
{
\frac{1}{\Delta} = \sum_{p\geq0} (-1)^p\frac{(2p!)}{4^p(p!)^2} A^{-(1/2+p)} V^p,
}\ee
a direct computation gives us the desired expression, truncated to a given degree in  $(\bx,\bxt)$.
As a result, the perturbative part of the Hamiltonian takes the form of a  power series of $(\bx,\bxt)$ whose monomials are expressed as:
\be
\begin{aligned}
 &\Gamma_{p_1,p_2,\pt_1,\pt_2,k_2}(\zeta_1,s_0) x_1^{p_1}x_2^{p_2}\xt_1^{\pt_1}\xt_2^{\pt_2} \exp(ik_2\zeta_2) \\
& \text{with} \quad  (p_1,p_2,\pt_1,\pt_2,k_2) \in \N^4\times\Z.
 \end{aligned}
 \label{eq:monom}
\ee
 Thanks to the symmetry of the problem relative to the invariant plane, certain coefficients $\Gamma_{p_1,p_2,\pt_1,\pt_2,k_2}$ cancel out.  Indeed, this symmetry requires that the parity of $p_1+p_2+\pt_1+\pt_2$  is the same as that of  $k_2$. The reader is referred to \cite{MaRoLa2002} for the proof of this property and more information about the reduction of the nodes.


\subsubsection{The average Hamiltonian}

A classical way to reduce the problem is to average the Hamiltonian over the fast angle $\zeta_2$. Resonant normal form theory \cite[see][]{{MeHa1992}} shows that there exists a canonical  transformation, close to the identity, that maps the Hamiltonian 
\be
\rH(\bZ, \bzeta, \bx, \bxt, \cC(s_0)) = \rH_K(\bZ) + \eps \rH_P(\bzero, \bzeta, \bx, \bxt, \cC(s_0)) 
\ee
to 
\be
\begin{aligned}
& F(\bZ, \bzeta, \bx, \bxt, s_0) = F_K(\bZ) + \eps F_P(\zeta_1, \bx, \bxt, s_0) +\gO(\eps^{3/2}) \\
 \qtext{with}   & \rH_K = F_K \qtext{and} 
 F_P(\zeta_1,\bx,\bxt,s_0) = \frac{1}{2\pi}\int_0^{2\pi} \rH_P(\bzeta,\bx,\bxt,s_0)  d\zeta_2 \, ,
\end{aligned}
\label{eq:H_moy}
\ee
where the same notations are used for the variables before and after averaging.  

Throughout the rest of Section \ref{sec:moyen}, we will assume that the planetary masses (or $\eps$) are small enough that we can neglect the remainder of order $\eps^{3/2}$ in (\ref{eq:H_moy}), which will amount to restrict ourselves to the Hamiltonian: 
\be
 \Fb(\bZ,\zeta_1,\bx,\bxt,s_0) = F_K(\bZ) + \eps F_P(\zeta_1,\bx,\bxt,s_0).
\label{eq:H_moy_rest} 
\ee
As mentioned in the last paragraph of \ref{sec:expansion}, since in the expression (\ref{eq:monom}), $k_2=0$, $F_P$ contains only terms of even total degree in $(\bx,\bxt)$.
As a consequence, the Hamiltonian (\ref{eq:H_moy_rest}) can be expanded in power series of $(\bx,\bxt)$ under the form: 
\be
\Fb(\bZ,\zeta_1,\bx,\bxt,s_0) = F_K(\bZ) + \eps\sum_{p\geq 0} F_P^{(2p)}(\zeta_1,\bx,\bxt,s_0),
\ee
where the $F_P^{(2p)}$ are polynomial of degree $2p$ in $(\bx,\bxt)$ whose coefficients depend on $\zeta_1$ and $s_0$.

\subsection{The circular inclined case}
\label{sec:circ}

\subsubsection{Integrable approximation}
 \label{sec:circ_integrable}

The average Hamiltonian $\Fb$ being even in the variables $(\bx,\bxt)$, the manifold $x_1 = x_2 = \xt_1 = \xt_2 = 0$  is invariant by the flow of this average Hamiltonian   \citep[see][for more details]{RoPo2013}.  On this manifold, the average orbits of the two co-orbital bodies are circular with non-zero mutual inclination. We consider, in this section, the restriction of the average Hamiltonian to this manifold, that is:
\be
F_K(\bZ) + \eps F_P^{(0)}(\zeta_1,\bzero,\bzero,s_0).
\ee  
In the average problem (averaged over the angle $\zeta_2$), the action $Z_2$ is constant. Since, according to (\ref{eq:H_Kep_Z}), the only effect of $Z_2$ is to translate the position of any equilibrium points, we can get rid of this variable without loss of generality by imposing $Z_2=0$. As a consequence, we are left with the study of the one parameter family of integrable Hamiltonian
\be
\begin{aligned}
F_c(Z_1,\zeta_1,s_0) &=  -\alpha Z_1^2 + \eps\beta  F_1(\zeta_1,s_0) \\
 &= -\alpha Z_1^2 + \eps\beta \left(
\left(1 - s_0^2\right) \cos\zeta_1 - \frac{1}{\pi}\int_0^{\pi} \frac{d\zeta_2}{\sqrt{ U_{s_0}(\bzeta)} }
\right), \\
  \text{with} \quad  U_{s_0}(\bzeta) &= 2 - 2\cos\zeta_1 + 2s_0^2 \left(\cos\zeta_1 - \cos(\zeta_1 + 2\zeta_2) \right), \\
  \alpha &= \frac{3{\ns}^{4/3}}{2\mu_0^{2/3}}\frac{m_1+m_2}{m_1m_2}  \qtext{and} \beta = {\ns}^{2/3}\mu_0^{2/3} \frac{m_1m_2}{m_0} . 
 \end{aligned}
 \label{eq:Ham_circ}
 \ee

Note that the integral appearing in expression (\ref{eq:Ham_circ}) can be expressed using the incomplete elliptic integral of first kind as:
\be
\begin{aligned}
\frac{1}{\pi}\int_0^{\pi} \frac{d\zeta_2}{\sqrt{ U_{s_0}(\bzeta)}} &= \frac{1}{2\pi \sqrt u}
\left[
  \cF\left(\frac{\zeta_1 + \pi}{2}, \frac{s_0^2}{\sqrt{u}}\right) 
- \cF\left(\frac{\zeta_1-\pi}{2}, \frac{s_0^2}{\sqrt{u}}\right)
\right] \\
\qtext{where} p &= \sin(\zeta_1/2), \, u = p^2 +s_0^2 -p^2s_0^2, \\ 
\qtext{and} \cF(x,y) &= \int_0^x \frac{dz}{\sqrt{1 - y^2\sin^2z}}.
\end{aligned}
\ee
The proof is provided in the appendix \ref{sec:appendix_Hamcirc}.

\subsubsection{The \VFLtitre in the average problem}
\label{sec:vfl_avg}
For a given $s_0 = \sin(J_0/2)$, the fixed points of the Hamiltonian (\ref{eq:Ham_circ}) satisfy the implicit equation 
\be
(Z_1,\zeta_1) = (0,\zeta_0) \qtext{with} \dron {F_1}{\zeta_1}(\zeta_0,s_0)=0.
\label{eq:equil_circ}
\ee
It is important to note that this equation is independent of masses. The same will be true for the position of the fixed points solutions of this equation, provided that the terms in $\eps^{3/2}$ are negligible. We will see in section \ref{sec:num} that when this assumption is not fulfilled, the shape of the family of fixed points, or rather periodic orbits, evolves as the “planetary” masses increase.     

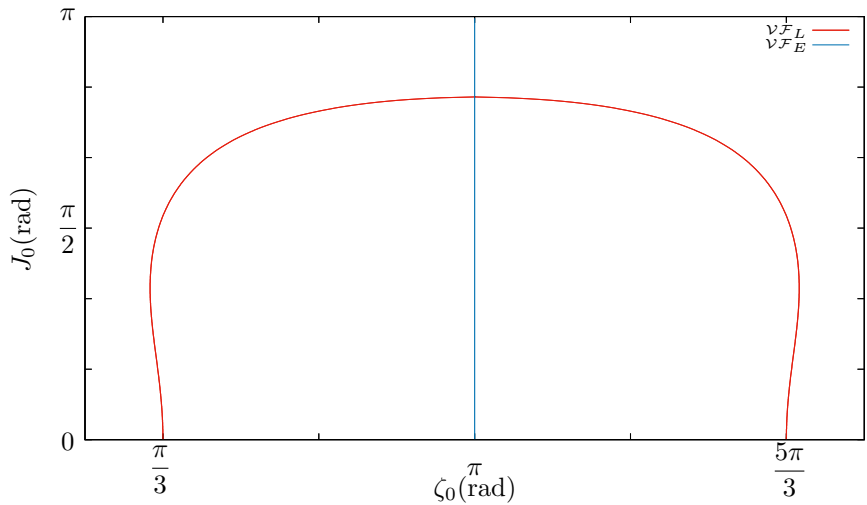
\begin{figure}[H]
\begin{center}
\begin{tikzpicture}
\newcommand{\xsca}{.9}
\newcommand{\ysca}{.75}
\newcommand{\dy}{-.14*\ysca}
\newcommand{\dx}{.15*\xsca}
\newcommand{\dxlab}{10.7*\xsca}
\node[xscale=\xsca, yscale=\ysca] at (0,0)  [above right] {\includegraphics[height=8cm,width=12cm]{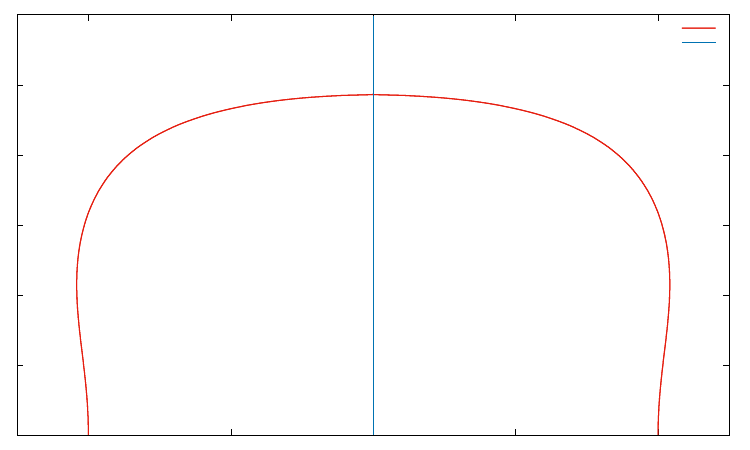} };
\path[black] (6.1*\xsca,-.5*\ysca) node  { $\zeta_0 (\rad)$  }; 
\path[black] (1.5*\xsca,\dy) node  { ${\displaystyle \frac{\pi}{3}}$ }; 
\path[black] (6.1*\xsca,\dy) node  { ${\displaystyle \pi}$ }; 
\path[black] (10.7*\xsca,\dy) node  { ${\displaystyle \frac{5\pi}{3}}$ }; 
\path[black] (-.5*\xsca,4.1*\ysca) node [rotate=90] { $J_0 (\rad)$ }; 
\path[black] (\dx,.35*\ysca) node  { $0$ }; 
\path[black] (\dx,4.1*\ysca) node  { ${\displaystyle \frac{\pi}{2}}$ }; 
\path[black] (\dx,7.9*\ysca) node  { ${\displaystyle \pi}$ }; 
\path[black] (\dxlab,7.6*\ysca) node  { {\tiny \VFL} }; 
\path[black] (\dxlab,7.35*\ysca) node  { {\tiny  \VFE} }; 
\end{tikzpicture}
\end{center}
\caption{The two families of fixed points derived from Lagrange equilibria in the average problem. The two branches of \VFL are identical up to permutation of the two small bodies.}
\label{fig:Marchal}
\end{figure}

Although, except for $s_0 = 0$, we cannot give an explicit expression for the solutions to this equation, we can either calculate them numerically or give an analytical approximation. In the latter case, we can express $\zeta_0$  in Taylor series of $s_0$, provided that the inclination is close enough to zero.
When $J_0 = s_0 = 0$, the fixed points correspond to Lagrange's equilateral equilibria at $\zeta_0 = \pi/3$ and $5\pi/3$ and to Euler colinear configuration at $\zeta_0 = \pi$. These equilibria are the starting point of fixed points families parameterized by the mutual inclination $J_0$ \changed{-- although we could also parametrize them with the total angular momentum $\cC$, $J_0$ provides a more natural and intuitive parameter}.
Solving equation (\ref{eq:equil_circ}) numerically, with initial condition $\zeta = \pi/3$ or $5\pi/3$ for  $J_0 = 0$, provides a unique family of fixed points, which will be denoted \VFL for Vertical Family of Lagrange. This family connects the two equilateral Lagrange plane configurations while containing inclined trajectories with maximum mutual inclination $J_c \approx 2.542567 \rad \approx  145.678^\circ$ at $\zeta = \pi$.
Since $F_1$ is even and $2\pi$-periodic,  \VFL is symmetric with respect to the line $\zeta_0 = \pi$. This line coincides with the second family of solutions of (\ref{eq:equil_circ}), which will be called the Vertical Family of Euler (\VFE). This latter structure originates from Euler's aligned configuration, in which the two planets are on both sides of the primary body, intersect \VFL at $J_c$ and persist until they meet at $J=\pi$, where the two planets encounter each other.  These two families, obtained by numerical resolution of the equation (\ref{eq:equil_circ}), are displayed in the figure \ref{fig:Marchal}.
This behavior of the families \VFL and \VFE has already been noticed in the average restricted 3-body problem by \cite{Marchal2009} and mentioned by \cite{Leleu2017PhD} in the case of the planetary 3-body problem.   According to \cite{Marchal2009}, the maximal value of the mutual inclination reaches on these families is equal to $145.68^\circ$ which is what we get here.
As mentioned above, although Marchal obtains this result in the average restricted problem, the agreement observed is due to the fact that the solutions to equation (\ref{eq:equil_circ}) are independent of the masses. This remains valid as long as the smallness of the planetary masses allows the terms of order $\eps^{3/2}$ in (\ref{eq:H_moy}).

The solutions of the equation (\ref{eq:equil_circ}) can also be approximation by Taylor expansion. 
Indeed \VFL can be expressed with a relative error lower than $10^{-5}$ up to  $J_0 \leq 60^\circ$ by the twentieth degree Taylor expansion (all analytical results were found using the TRIP language \cite{gastineauTRIPComputerAlgebra2011}):
\be
\begin{aligned}
\zeta_c(J_0) &= \frac{\pi}{3} -  \sqrt3 \left( 
\frac{1}{3} s_0^2
-  \frac{4}{9} s_0^{4} 
 - \frac{17}{27} s_0^{6}
 + \frac{1865}{2592} s_0^{8}
 + \frac{13877}{9720} s_0^{10}
 - \frac{75553}{23328} s_0^{12} \right. \\
 &\left. - \frac{63599}{10206} s_0^{14}
 + \frac{254313827}{17915904} s_0^{16} 
 + \frac{384292411}{13436928} s_0^{18}
 - \frac{86371296841}{1209323520} s_0^{20}
\right) \\
 & + \gO(s_0^{22}),
 \end{aligned}
 \label{eq:zeta_taylor}
\ee
where $s_0 = \sin(J_0/2)$. 
The expression of the branch originating from $L_5$ is expressed in the same way by taking
$\zeta'_c(J_0) = 2 \pi - \zeta_c(J_0)$.

Note that by expressing $\zeta_c(J_0)$  using directly $J_0$ rather than $s_0$, we find 
\be
\zeta_c(J_0) \approx  \frac{\pi}{3} -   \frac{\sqrt3}{12}J_0^2  + \frac{5\sqrt3}{144}J_0^4 + \cdots
\ee
which corresponds to the expression given by \cite{NaMu2000} for the restricted 3-body problem.

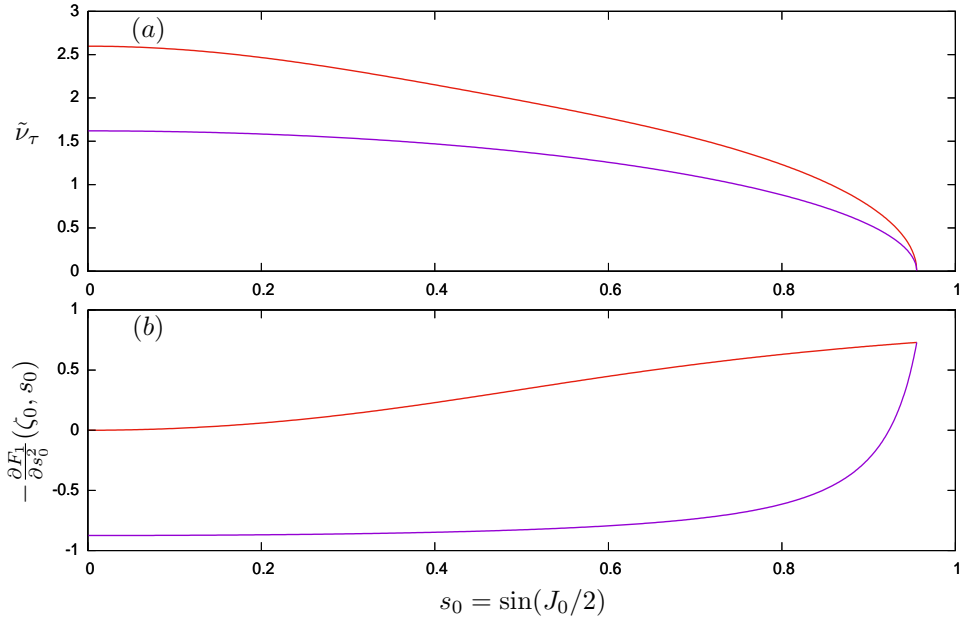
\begin{figure}[H]
\begin{center}
\begin{tikzpicture}
\newcommand{\xsca}{.88}
\newcommand{\ysca}{.88}
\newcommand{\dy}{7*\ysca}
\newcommand{\dx}{1.8*\xsca}
\newcommand{\dxlab}{10.7*\xsca}
\node[xscale=\xsca, yscale=\ysca] at (0,0)  [above right] {\includegraphics[height=9cm,width=14cm]{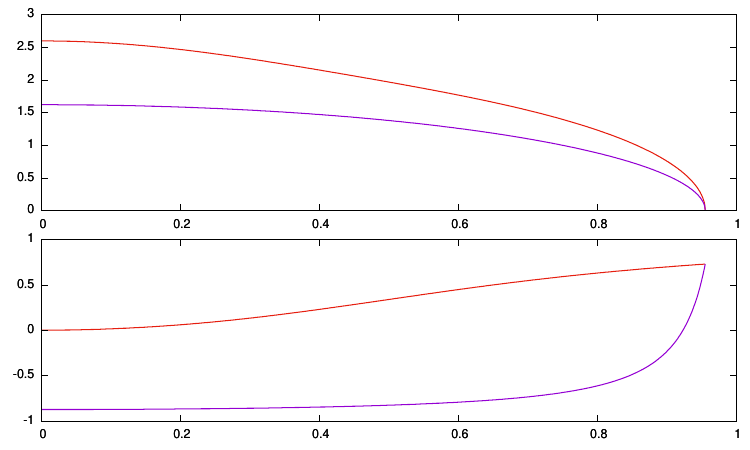} };
\path[black] (\xsca*7.4,-.08*\ysca)  node  { $s_0 = \sin(J_0/2)$ }; 
\path[black] (0*\xsca,\dy)  node  { $\nut_{\tau}$ }; 
\path[black] (\dx,8.55*\ysca)  node  { $(a)$ }; 
\path[black] (0*\xsca,2.5*\ysca)  node  [rotate=90] { $-\dron{F_1}{s_0^2}(\zeta_0,s_0)$ }; 
\path[black] (\dx,4.05*\ysca)  node  { $(b)$ }; 
\end{tikzpicture}
\end{center}
\caption{(a): Behavior of eigenvalues and precession frequencies of nodes along families \VFE (purple) and \VFL (red). (b): Precession of the ascending node. }
\label{fig:freq_FL34}
\end{figure}

As long as we limit ourselves to studying the dynamics restricted to the plane $x_j = 0$,
the linear stability of the fixed points is given by the eigenvalues of the linearized differential system at the equilibrium these are roots of the characteristic polynomial:
\be
P(X) = X^2 -  2 \eps\alpha\beta \frac{\partial^2F_1}{\partial \zeta^2} (\zeta_0)
= X^2 - 3\eps{\ns}^2\,\frac{m_1 +m_2}{m_0}\frac{\partial^2F_1}{\partial \zeta^2} (\zeta_0).
\label{eq:stab_circ}
\ee
The nature of the fixed point (elliptic/hyperbolic) is deduced from the sign of  $\frac{\partial^2F}{\partial \zeta^2} (\zeta_0)$.
Below the critical mutual inclination $J_c$, equilibria belonging to the families \VFL are stable (purely imaginary eigenvalues), whiles these are unstable (real eigenvalues) along this part of \VFE. These fixed points degenerate (zero eigenvalues) at the bifurcation point where the three families join together \changed{(\VFE and the two branches of \VFL, symmetrical under permutation of the two small bodies)}. 
Above this critical value,  \VFE becomes elliptic and ends  with a collision orbit at $J_0 = 180^\circ$.  
The eigenvalues of the elliptic family \VFL, resp. hyperbolic \VFE  of the form $\pm i\nu_L$, resp. $\pm \nu_E$, with 
\be
\nu_\tau =  \sqrt{\eps}\ns\nut_\tau\sqrt{\frac{m_1+m_2}{m_0}} \qtext{and} 
\nut_\tau = \sqrt{\gamma_\tau}\left( 1 + \gO(s_0^2)\right)\qtext{and} \tau \in \{L,E\}
\label{eq:fre_nu}
\ee
where $\gamma_L = 27/4$ and $\gamma_E = 21/8$. 
The dimensionless quantities $\nut_\tau$, for which Taylor expansions are given in the appendix \ref{sec:appendix_Taylor_freq}, are plotted  in figure \ref{fig:freq_FL34}.a.
The range of  $J_0$ is limited to the interval $[0,J_c]$ as higher values are of limited relevance for the present study. $\nut_L$ (red)  and $\nut_E$ (purple) both decrease until they reach zero at $J_c$, where the three families merge. 
The Jacobi's reduction allowed us to eliminate the inclination and longitude of the node of each planet. As indicated in section \ref{sec:Jacobi}, theses quantities can be derived from the expressions (\ref{eq:inclin})  and (\ref{eq:noeud})  (to within a constant for $\Omega_p$). 
Thus, the ascending node of the orbits of \VFL precess at a constant rate, which, according to (\ref{eq:equil_circ}), is given by:
\be
\dot \Omega_p =\dot \psi_p =  \dron{F_c}{\cC}(0,\zeta_0,s_0) = -\eps\ns\frac{m_1+m_2}{2m_0}
\sqrt{1 - 4s_0^2\frac{m_1m_2}{(m_1+m_2)^2}} \dron{F_1}{s_0^2}(\zeta_0,s_0).
\label{eq:precession_noeud}
\ee
The function $-\dron{F_1}{s_0^2}(\zeta_0,s_0)$ is plotted against $s_0 = \sin(J_0/2)$ in the figure \ref{fig:freq_FL34}.b and its Taylor expansion is provided in the appendix \ref{sec:appendix_Taylor_freq}. 
On the family derived from the equilateral configuration, the node line has a prograde precession motion, whereas  on \VFE, the motion is retrograde up to $s_0\approx .9235$  (or $J_0\approx 135^\circ$), then prograde for higher values.

\subsubsection{Comparison with numerical simulations in the non-average problem}
\label{sec:compar}

\begin{figure}[H]
\begin{center}
\begin{tikzpicture}
\newcommand{\xsca}{.88}
\newcommand{\ysca}{.85}
\newcommand{\dy}{7*\ysca}
\newcommand{\dx}{1.8*\xsca}
\newcommand{\dxlab}{10.7*\xsca}
\node[xscale=\xsca, yscale=\ysca] at (0,0)  [above right] {\includegraphics[height=9cm,width=14cm]{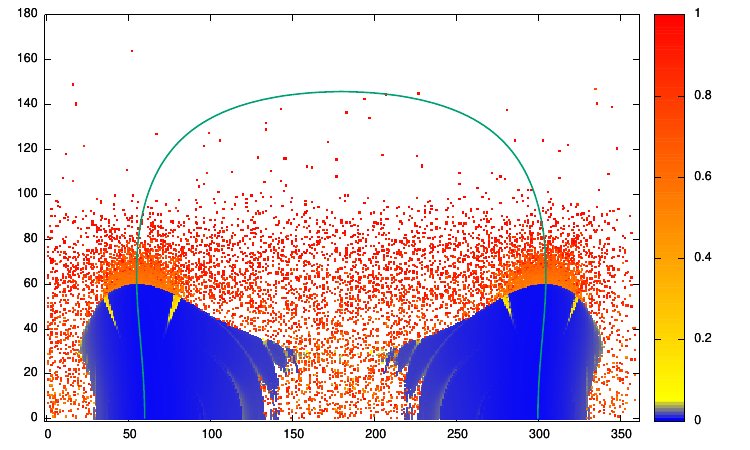} };
\path[black] (\xsca*6.5,-.08*\ysca)  node  { $\lam_1 - \lam_2$ (deg)}; 
\path[black] (0*\xsca,4.8*\ysca) node  [rotate=90] { $J_0$ (deg) }; 
\path[black] (13.7*\xsca,4.8*\ysca) node  [rotate=90] { max($e_1,e_2$) }; 
\end{tikzpicture}
\end{center}
\caption{Stability map of the co-orbital region in the plane ($\lambda_1 - \lambda_2, J_0)$ for $m_1 = m_2 = 10^{-3}m_0$. The stability index, from blue to red, corresponds to the maximum eccentricity values reached by the two planets during their integration in the full problem, while the white is associated with orbits that either left co-orbital resonance before the end of integration or collided. The green curve is the projection of \VFL (computed in the average problem) onto the stability map.}
\label{fig:Moria}
\end{figure}
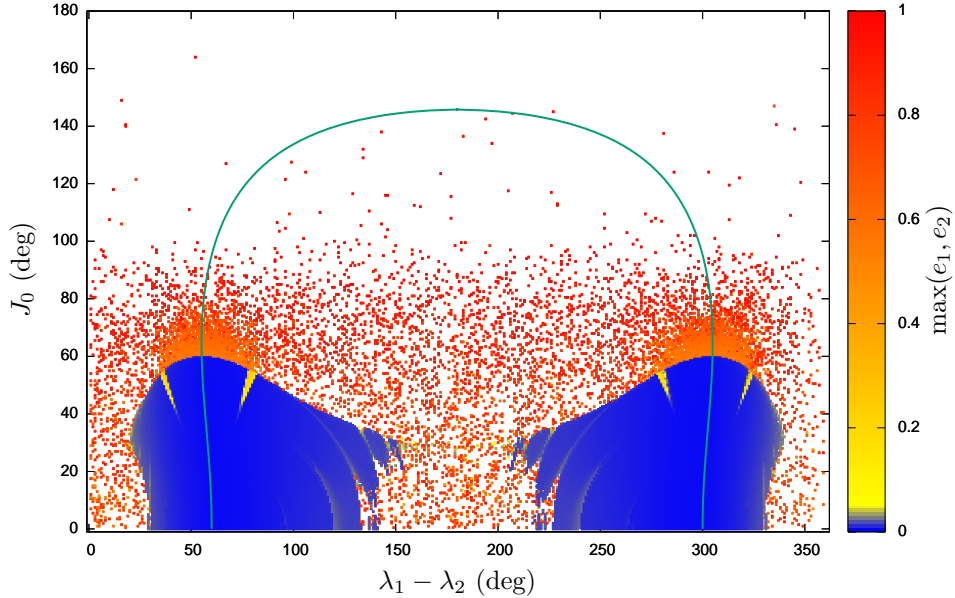

Let us now see how this local aspect of the geometry of the problem can be part of a global view of its dynamics, and in particular how the presence of \VFL can influence the stability of co-orbital regions. First, let us examine the stability map shown in the figure \ref{fig:Moria}.
This map, consists of a grid of initial conditions with coordinates $(\lam_1 - \lam_2,J_0)$ with $\lam_1 \in [0^\circ,360^\circ]$ and $2I_1 = 2I_2 = J_0 \in [0^\circ,180^\circ]$, the other initial conditions and parameters being set as follows: $a_1 = a_2  = 1$, $e_1 = e_2=0$, $\lam_2=0$, $\Om_1 = \Om_2 + \pi = 0$. and $m_1/m_0 = m_2/m_0 = 10^{-3}$.
 Each point on the grid is associated with a color indicating the maximum value reached by $e_1$ and $e_2$ during a numerical integration lasting $40406$ orbital revolutions. As indicated on the color bar on the right-hand side, blue corresponds to maximum values below $0.05$, while colors ranging from yellow to dark red correspond to values between $0.05$ and $1$. White is associated with systems that either left co-orbital resonance before the end of integration or collided. 
 Here we use the maximum eccentricity as an indicator of instability. On the one hand, this is a relative indicator in the sense that it is the variations in its value that indicate changes in dynamic behavior such as the presence of resonance or chaotic regions \citep[see][]{Leleu2017PhD, ZhDvSu2009}. On the other hand, it also provides absolute information in the sense that it indicates the regions in which circular or quasi-circular orbits exist.
The contrast between regions with regular dynamics (blue) and chaotic dynamics is very marked. 
The orbits that are located at the center of the two co-orbital regions ($\zeta$ close to $60^\circ$ or $300^\circ$)  undergo only very slight variations in eccentricity (smaller than $0.05$), and over the time interval considered, are indistinguishable from quasi-periodic orbits; \changed{the measurement of the diffusion in the frequency space provided by frequency map analysis \cite{Laskar1990, Laskar1999} allows us to confirm this.}
Moving horizontally (variation of $\zeta$) or vertically (increasing $J_0$) on this map, we encounter a very narrow and sharp transition zone between stable and chaotic regions. Whether, in the horizontal direction, this boundary is caused by the overlap of secondary resonances associated with commensurabilities between the libration frequency (associated with the resonant angle $\lam_1 -\lam_2$) and the orbital frequency common to both planets \citep[see][]{ErNaSaSuFro2007,RoBo2009}, the reason for the presence of the vertical transition zone at 60 degrees in the figure \ref{fig:Moria} is not clearly identified. 

Before addressing this question, let us return to our \VFL family.
Since the members of this family are fixed points of the reduced average problem, they should correspond, in the reduced full problem (non-average and moving with the line of nodes), to periodic orbits whose frequency is equal to the mean motion common to the two planets. In a fixed reference frame, these are therefore quasi-periodic orbits whose second frequency is the precession frequency of the nodes. Let us add that the elliptic elements, which are constant in the average problem, are animated by a fast (mean motion) motion whose amplitude is of the order of the perturbation $\eps$. In particular, as we will see in figure \ref{fig:coords_along_traj}, the trajectories may no longer be circular, but their eccentricity also remains of the order of $\eps$. For these reasons, although they are not expressed in the same coordinate system, it remains relevant to project \VFL (green curve) onto the stability map displayed on \ref{fig:Moria}. As long as the mutual inclination along \VFL remains below $60^\circ$, the family passes through the core of resonance (most stable regions in blue); above this value, all \VFL trajectories become unstable\footnote{A similar observation is made for elliptical orbits in \cite{GiLe2016}, although the cause of the instability is not identified.}\!\!.

One of the keys to this puzzle lies in calculating the eigenvalues of the equilibrium points of the family. We briefly mentioned this issue in Section \ref{sec:vfl_avg}, but the eigenvalues calculated there only corresponded to those associated to the dynamics in the plane $(\zeta,Z)$. We must therefore consider the dynamics in the $(\bx,\bxt)$ planes of coordinate; this would require numerically solving the variational equations for variations in the $(\bx,\bxt)$ directions.
%
Instead, in the following section we will directly numerically compute the periodic orbits (in rotating coordinates) in the full problem and solve their variational equations. 
The characteristic exponents thus calculated will provide information about the stability of the orbits.

\changed{Moreover, since we are only considering the average Hamiltonian at order $1$ of the masses, the equation (\ref{eq:equil_circ}) is independent of the masses. Of course, this approximation is only valid in the context of planetary problems (small masses). Determining the upper limit of planetary masses for which this approximation is valid would require at least calculating the average Hamiltonian to the second order, which is beyond the scope of this paper. Direct numerical integration of the full, non-approximated problem will circumvent this limit, and give us some insight on the validity of the approximation: Figure \ref{fig:comparison_avg_full} compares between results in the average and the full problem for various mass parameters (although the angles used in the two problems are slightly different, section \ref{vfl_planetary} for a more detailed discussion).}

To answer both of these problems, in the following section we will directly numerically compute the periodic orbits (in rotating coordinates) in the full problem and solve their variational equations. 
The characteristic exponents thus calculated will provide information about the stability of the orbits.

\section{Numerical search of \VFLtitre in the full problem}
\label{sec:num}

\subsection{Reformulating the problem: the Hill variables}

\subsubsection{Introducing the Hill variables}
\label{sec:hill_coords}

Our goal is now to show the existence of the \VFL family in the full — i.e., non-average problem and to study the stability of trajectories along this family. However, in the Poincaré coordinates used before, the Hamiltonian and its associated vector field cannot be explicitly expressed without truncating in eccentricity.  To avoid this limitation, we must introduce coordinates that allow an explicit expression of the Hamiltonian and are adapted to the Jacobi reduction. Indeed, as mentioned in section \ref{sec:compar}, in the draconic reference frame, the orbits sought are not quasi-periodic but simply periodic, which will allow us to use classical methods for finding periodic orbits.

The Hill variables \citep[see][and references therein]{Laskar2017Andoyer}, given by $(r_j, w_j, \Omega_j, R_j, G_j, G \cos{I_j})_{j = 1, 2}$ satisfy this double constraint. The draconic true anomalies $w_j$, longitudes of the nodes $\Om_j$, angular momenta $G_j$ and vertical angular momenta $G_j \cos{I_j}$ ($I_j$ being the inclinations) of each planets have already been introduced in \ref{sec:Jacobi}, while $r_j$ and $R_j$ are respectively the distance between the bodies of masses $m_j$ and $m_0$, and  the norm of theirs linear momentum in the barycentric frame.
Similarly to \ref{sec:Jacobi}, we perform the Jacobi reduction and eliminate the nodes and vertical momenta, but keep the angular momentum $C$ as a parameter. As a consequence,  $8$ dimensions remain corresponding to the coordinates: $(r_j, w_j, R_j, G_j)_{j = 1, 2}$ to which the parameter $\cC$ has to be added.
We can express the positions and momenta in the canonical heliocentric coordinates as:

\begin{equation}
    \br_j = R_x(\Om_j) R_z(I_j) R_x(w_j) \begin{pmatrix}
      r_j \\ 0 \\ 0
    \end{pmatrix} , \quad
    \brt_j = R_x(\Om_j) R_z(I_j) R_x(w_j) \begin{pmatrix}
      R_j \\ G_j / r_j \\ 0
    \end{pmatrix}.
  \label{eq:hill_expr}
\end{equation}
where for $\theta \in \R$, $R_u(\theta)$ is the rotation matrix of angle $\theta$ along axis $u$:
\begin{equation*}
    R_x(\theta) = \begin{pmatrix}
        1 & 0 & 0 \\
        0 & \cos{\theta} & -\sin{\theta} \\
        0 & \sin{\theta} & \cos{\theta}
    \end{pmatrix} ,
    R_y(\theta) = \begin{pmatrix}
        \cos{\theta} & 0 & \sin{\theta} \\
        0 & 1 & 0 \\
        -\sin{\theta} & 0 & \cos{\theta}
    \end{pmatrix} ,
    R_z(\theta) = \begin{pmatrix}
        \cos{\theta} & -\sin{\theta} & 0 \\
        \sin{\theta} & \cos{\theta} & 0 \\
        0 & 0 & 1
    \end{pmatrix}
\end{equation*}

The planetary Hamiltonian (\ref{eq:ham_cart}) can now be written in Hill variables. As $\brt_j^2 = R_j^2 + G_j^2/r_j^2$ and $\br_j^2 = r_j^2$, its Keplerian part $H_K$ is trivially translated in Hill coordinates. In order to express the perturbative part in terms of Hill variables, it suffices to note that the dot product $\brt_1 \cdot \brt_2$ can be written as:
\begin{align*}
    \brt_1 \cdot \brt_2 = \left( R_x(w_1) \begin{pmatrix}
        R_1 \\ G_1 / r_1 \\ 0
    \end{pmatrix} \right) \cdot \left(
        R_x(-I_1) R_z(\Om_2 - \Om_1) R_x(I_2) R_x(w_2)
        \begin{pmatrix}
            R_2 \\ G_2 / r_2 \\ 0
    \end{pmatrix} \right),
\end{align*}

As $\Om_2 - \Om_1 = \pi$, the mutual inclination $J$ is the sum of the two individual inclinations. Thus, a final straightforward calculation establishes the relation $$R_x(-I_1) R_z(\pi) R_x(I_2) = R_z(\pi) R_x (J),$$
which leads to the final expression:
\be
\begin{aligned}
\brt_1 \cdot \brt_2 = &  \left(\cos{w_1} \rt_1 + \sin{w_1} \frac{G_1}{r_1}\right) \left(  -\cos{w_2} \rt_2 + \sin{w_2} \frac{G_2}{r_2} \right) \\
                      & + \cos J\left( \sin{w_1} \rt_1 - \cos{w_1} \frac{G_1}{r_1}\right) \left( \cos{w_2} \frac{G_2}{r_2} - \sin{w_2} \rt_2\right).\\
\end{aligned}
\ee
Using the same transformation for the dot product $\br_1 \cdot \br_2$, the planetary Hamiltonian is expressed in Hill variables as follows:
\begin{multline}
  H_\cC(r_j, w_j, \Omega_j, R_j, G_j) =
  \sum_{i=1}^2 \left(
    \frac{m_0 + \eps m_i}{2 m_0 m_i}\left( R_i^2 + \frac{G_j^2}{r_j^2}\right)- \cG \frac{m_0 m_i}{r_i}
  \right) \\
  + \frac{\eps}{m_0} \left[
    - \cos{w_1}\cos{w_2} R_1 R_2
    - \sin{w_1}\sin{w_2} \frac{G_1 G_2}{r_1 r_2} \right. \\
    + \cos{w_1}\sin{w_2} \frac{G_2 R_1}{r_2}
    + \sin{w_1}\cos{w_2} \frac{G_1 R_2}{r_1} \\
    + \cos J \left(
      -\sin{w_1}\sin{w_2} R_1 R_2
      - \cos{w_1}\cos{w_2} \frac{G_1 G_2}{r_1 r_2} \right. \\ \left. \left.
      - \sin{w_1}\cos{w_2} \frac{G_2 R_1}{r_2}
      - \cos{w_1}\sin{w_2} \frac{G_1 R_2}{r_1}
  \right)
  \right] \\
  - \eps \cG m_1 m_2 \left(
    r_1^2 + r_2^2 - 2 r_1 r_2 \cos{w_1}\cos{w_2}
    - 2 r_1 r_2 \sin{w_1}\sin{w_2} \cos J
  \right)^{-1/2},
  \label{eq:ham_hill}
\end{multline}
with $$\cos J =  \frac{C^2 - G_1^2 - G_2^2}{2 G_1 G_2}. $$
Note that, even though the longitudes of the nodes do not appear explicitly in the expression of the Hamiltonian, we can recover the node precession (the two node precessions are equal, as $\Om_1 = \Om_2 + \pi$) by integrating (\ref{eq:noeud}) along a trajectory $\mathbf{x}(t)$ by:
\begin{equation}
    \forall t, \Om (t) = \Om(0) + \int_0^t
        \frac{\partial }{\partial \cC}\big(H_\cC (\mathbf{x}(\tau))\big) d\tau.
\end{equation}


\paragraph{Limit cases}

The Jacobi reduction has two limit cases: the planar case ($J=0$) and the restricted case ($m_j = 0$ for one $j$).

In the planar case, the nodal line -- the intersection of the orbital plane and the invariant plane -- is ill-defined, as the planes are identical. However, no singularity appears in the Hamiltonian, nor in the associated vector field; therefore, we can still use the Hill variables in this case. Let us remember that in this case the only relevant angle is the true longitude $\Om_j + w_j$.

The restricted limit, where one of the masses is zero, is not as smooth: the Jacobi reduction simply cannot be used there, as the zero-mass body isn't constrained by the direction of the total angular momentum. Thus, we restrict ourselves to low, but non-zero masses.

\subsubsection{The Lagrange relative equilibrium}

 As indicated above, Jacobi reduction still allows us to consider coplanar motions and, in particular, Lagrange's relative equilibria. Thus, in Hill coordinates, the quantities $(r_j,R_j,G_j,w_1-w_2)$ are constants given by: 
\be
\begin{aligned}
    r_1 &= \rho,                                                   &r_2 &= \rho, \\
    R_1 &= -\frac{\sqrt{3} \om \rho \eps m_1 m_2}{2(m_0 + m_1 + m_2)},  &R_2 &= \frac{\sqrt{3} \om \rho \eps m_1 m_2}{2(m_0 + m_1 + m_2)}, \\
    G_1 &= \frac{\om \rho^2 m_1 (2m_0 + m_2)}{2(m_0 + m_1 + m_2)},  &G_2 &= \frac{\om \rho^2 m_2 (2m_0 + m_1)}{2(m_0 + m_1 + m_2)}, \\
    &w_1 - w_2 = \frac{5\pi}{3},
\end{aligned}
\label{eq:L4_Hill}
\ee
where $\rho$ is the radius of the circles described by the two planets, which can be considered a scale factor, and $\om$ is the planetary mean motion. These two quantities are, of course, not independent but linked by Kepler's third law. 
Contrary to their difference, the true anomalies $w_1, w_2$ are not constant. Indeed, they satisfy the first order differential equation:
\be
\dot w_j = \om\left(a + b\sin(2w_j+ (-1)^j\pi/3) + c\cos(2w_j+ (-1)^j\pi/3)\right),
\label{eq:edo_w}
\ee
where $a, b, c$ are constants depending only on the three masses (see appendix \ref{sec:appendix_Lagrange}) that satisfy the relation:  $a^2 - b^2 - c^2 = 1$. Thanks to this last relation, it can be shown that the solution to this differential equation takes the form: 
$w_j(t) = \omega t + F(\omega t) $, where $F$ is $2\pi$-periodic. Although the $w_j$ are pulsating (due to $F(\om t)$), the physical angles to consider in the planar case are the $\Om_j+w_j$, which can be shown to be exactly equal to $\om t$.
%

The spectral stability of the relative equilibrium described above can be deduced from the eigenvalues of the monodromy matrix of its variational equation. The corresponding spectrum is given by: 
\begin{equation}
    \begin{pmatrix}
        1 \\
        \exp{
        \left(
            \pm \frac{i \om}{\sqrt2} \sqrt{1 - \sqrt{1 - \frac{27 p}{\sigma^2}}}
        \right)} \\
        \exp{\left(
            \pm \frac{i \om}{\sqrt2} \sqrt{1 + \sqrt{1 - \frac{27 p}{\sigma^2}}}
        \right)}
    \end{pmatrix}
    \label{eq:lagrange_eigenvalues}
\end{equation}
where $\sigma = m_0 + m_1 + m_2$ and $p = m_0 m_1 + m_0 m_2 + m_1 m_2$. The eigenvalue $1$ is of multiplicity $4$ -- two because of the Hamiltonian's time invariance, and two because of the homographic family. We can also find the Gascheau criterion here: if $27 p / \sigma^2$ is greater than $1$, the root becomes complex, and the four non-$1$ eigenvalues aren't of modulus $1$ anymore; the Lagrange relative equilibrium goes from being elliptical to hyperbolic-elliptic and thus unstable.

\subsection{Methods: searching for periodic orbits}
\label{sec:num_methods}
As mentioned in section \ref{sec:compar}, after reduction of the nodes, the orbits that we are looking for are periodic. Due to the three-body problem's scale invariance, this period is arbitrary since resizing an orbit by a factor $\alpha$ corresponds to changing the time scale by a factor of $\alpha^{2/3}$. We will impose that this period $T$ be equal to $1$.

The periodic orbits that are members of the \VFL family can be found out as fixed points of a first return map.
Here, we consider the function $f$ defined as follows:
\begin{equation}
  f : \begin{cases}
      \R^9 \rightarrow \R^8 \\
      (\mathbf{x}; \cC) \mapsto \Phi_T^\cC(\mathbf{x}) - \mathbf{x}
  \end{cases}
\end{equation}
where $\mathbf{x} = (R_j, G_j, r_j, w_j)_{j=1,2}$ and $\Phi_T^\cC$ is the flow at time $T$ of the vector field associated to the Hamiltonian $H_\cC$ defined in (\ref{eq:ham_hill}). 
The Hamiltonian system possesses two independent first integrals: the energy and the total angular momentum. This diminishes the rank of $f$ by two, meaning the roots will never be locally unique. This can be solved by adding two sections: instead of finding roots of $f$, we search for roots of\footnote{Note that achieving local uniqueness can also be done via a change of coordinates: for instance, instead of \say{removing} one body in the coordinates, we could have added a six-dimensional section $\sigma: (\mathbf{p}, \mathbf{q}) \mapsto (\mathbf{p_0}, \mathbf{q_0})$ setting the first body at the origin. In this example, adding sections is of course computationally worse, but in some cases, we don't know how to define such a change of coordinates. For our case, we used the general properties (symmetries and invariances) of the three-body problem in changes of coordinates, and for matters specific to this family, we define sections.}
\begin{equation}
  (f, \sigma_1, \sigma_2)
  : \mathbb{R}^9 \rightarrow \mathbb{R}^{8+2}
  \label{eq:sigma1}
\end{equation}
The first section, for isolating a point on the orbit (also called a Poincaré section), is
\begin{equation}
  \sigma_1(\mathbf{x}) = w_1 + w_2
  \label{eq:sigma2}
\end{equation}
to keep the symmetry between the two small bodies.

For the family section, given the insight of \cite{Marchal2009} and \cite{Leleu2017PhD}, we choose the mutual inclination $J$\footnote{The parameter can be challenging to find, especially if there is a turning point. Techniques, such as pseudo-arc-length continuation, exist to circumvent such problems; however, we didn't need them for the results presented in this article.}, which can be expressed with the reduced Hill variables as
\begin{equation}
  \sigma_2(\mathbf{x}; \cC) = \cos J - \cos J_p = \frac{\cC^2 - G_1^2 - G_2^2}{2 G_1 G_2} - \cos{J_p}
\end{equation}
where $J_p$ is a parameter we will increase along the family. Note that, contrary to the average problem, $J$ is \textit{not} constant along the orbit; $J_p$ is the value of $J$ \textit{at the Poincaré section $\sigma_1$}. We will numerically see that this corresponds to the minimum of $J$ along the orbit (cf. Fig. \ref{fig:coords_along_traj}). For a summary of the different $J$ parameters used, see table \ref{tab:J_params}.

\begin{table}[htpb]
    \centering
    \begin{tabular}{c|c|c|c}
         \thead{\textbf{Name}} & \thead{\textbf{Case}} &
         \thead{\textbf{Definition}} & \thead{\textbf{Observation along}\\ \textbf{the \VFL family}} \\
         \hline
         $J_0$ & Average problem & $J$ at zero eccentricity for given $\cC$ & Maximum $J$ for given $\cC$ \\
         $J_p$ & Full problem & $J$ on an orbit when $w_1 + w_2 = 0$ & Minimum $J$ along the orbit
    \end{tabular}
    \caption{\changed{Summary of the $J$ parameters used in the article. $J_0$ is defined in section \ref{sec:expansion}, and $J_p$ in section \ref{sec:num_methods}.}}
    \label{tab:J_params}
\end{table}

We will focus on the case of two equal masses, $m_1 = m_2$, giving us the triplet $(1-2\eps, \eps, \eps)$ -- motivated by (\ref{eq:equil_circ}), showing that if the masses are small, the periodic orbits are independent of them; see \ref{sec:appendix_unequal} for a small discussion about unequal masses.

To compute $\Phi_T$, we use the Taylor method, via the \texttt{taylor} software from \cite{gimenoNumericalComputationHighorder2022}; this gives access to the numerical derivatives of the flow to an arbitrary order (computed with automatic differentiation). We embed the calls to \texttt{taylor} in a Julia \citep{bezanson2017julia} workflow, thanks to the \texttt{TaylorInterface.jl} (\url{https://github.com/alseidon/TaylorInterface.jl}) package, developed by the authors. \changed{\texttt{taylor} allows automatically chooses step size and approximation order when specified both a relative and absolute maximum error; we set both to $10^{-16}$.}

The root search can be done with a classical multi-dimensional Newton-Raphson algorithm; \changed{we choose a tolerance of $10^{-12}$. As we are searching a continuous family of orbits, we use the previous point as an initial guess for the next search (beginning with the Lagrange equilibrium), and allow for a maximum of $50$ iterations (otherwise, the search fails).}

The algorithm requires computing the function's first order derivatives, which are provided by the \texttt{taylor}; to compute the non-square Jacobian's inverse (with sections), we use QR decomposition. This also comes in handy for getting the stability of the orbit once found. Following \cite{MeHa1992}'s nomenclature, the Jacobian $\Phi_T$ is called the \textit{monodromy matrix}, and its eigenvalues are the \textit{characteristic exponents} of the orbit. They come as quadruplets $(\lambda, \Bar{\lambda}, \lambda^{-1}, \Bar{\lambda}^{-1})$ with $\lambda \in \mathbb{C}^*$. If all eigenvalues are of amplitude $1$, the orbit is elliptic and thus stable; otherwise, it is hyperbolic-elliptic, and unstable.

\subsection{Numerical results for \VFLtitre in the full problem}
\subsubsection{\VFLtitre for \texorpdfstring{$\eps = 10^{-3}$}{ε=1e-3}}
\label{vfl_planetary}

\begin{figure}[htbp]
  \centering
  \includegraphics[width=0.9\textwidth]{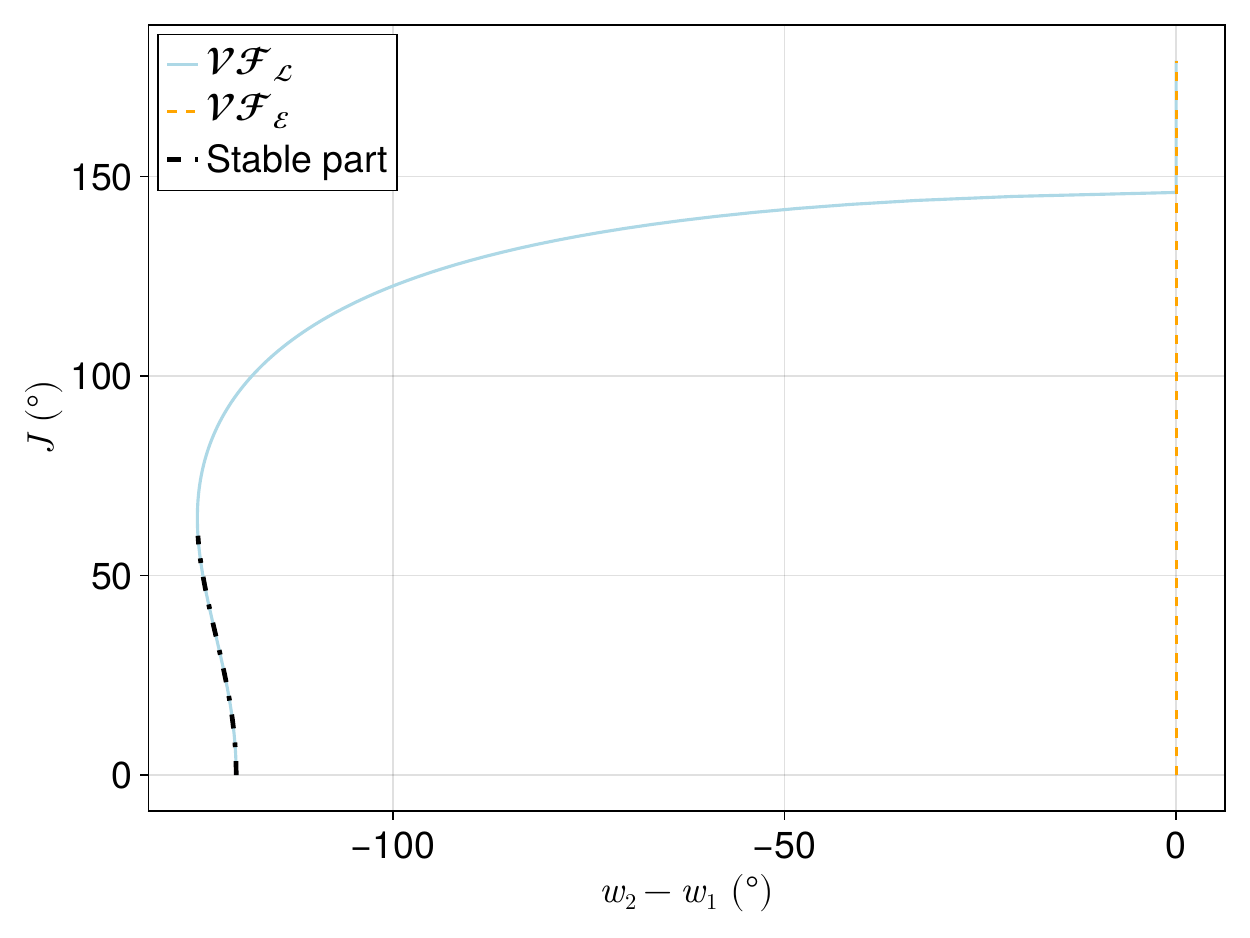}
  \caption{\VFL and \VFE for $\eps = 10^{-3}$, in the $( w_1 - w_2, J)$ plane.}
  \label{fig:marchal_small_masses}
\end{figure}

For $\eps = 10^{-3}$, we find the family for $J \in [0^\circ; 179^\circ]$; we conjecture these inclined co-orbitals lead to double collisions at $J = 180^\circ$ (the two small masses follow retrograde orbits of identical semi-major axis). Figure \ref{fig:marchal_small_masses} is the equivalent in the full problem of \citep[fig. 3]{Marchal2009} in the restricted problem and Fig. \ref{fig:Marchal} in the average problem; \changed{although the angle difference is different in each case: Fig. \ref{fig:marchal_small_masses} considers draconic true anomalies $w_j$; \citep[fig. 3]{Marchal2009} considers mean longitudes; and Fig. \ref{fig:Marchal} considers draconic mean anomalies $\lat_j$.
However, the variations between these angles are of order of the eccentricities, which are of order $\eps$, meaning the quantitative difference isn't visible; further comparison between $w$ and $\lat$ can be found in Appendix \ref{sec:diff_avg_full}.}

As Marchal noted in the restricted case, the family joins with a family emerging from the Euler relative equilibrium ($L_3$ in the restricted problem), in a bifurcation at about $J = 145^\circ$. This bifurcation is also connected to the Lagrange equilateral configuration  with permuted small masses ($L_4$ vs $L_5$), and continues up to the conjectured collisions.

\begin{figure}
    \centering
    \includegraphics[width=0.9\linewidth]{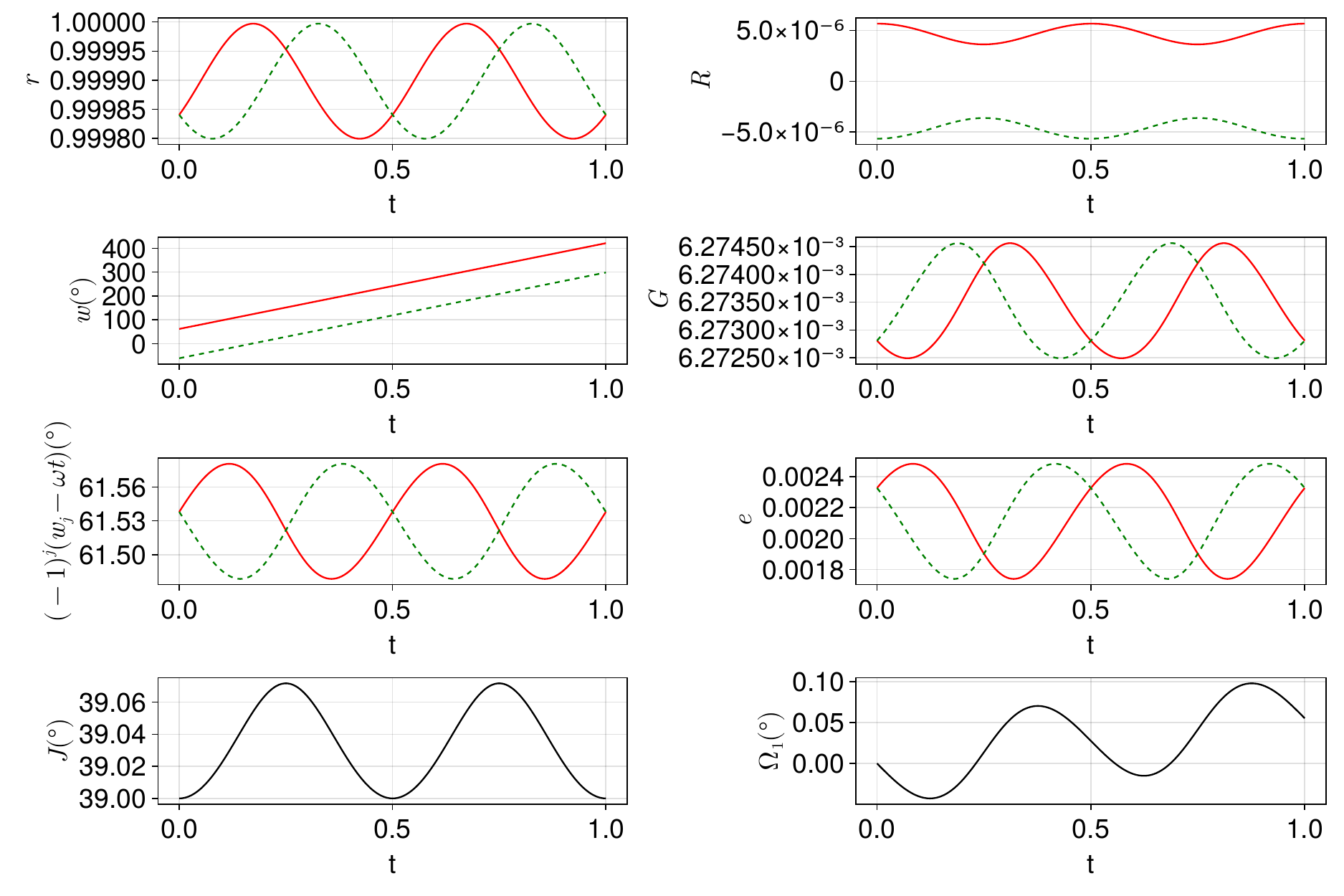}
    \caption{Coordinates for one trajectory, at $\eps = 10^{-3}$ and $J = 39^\circ$. We can see that the initial J, fixed as $J_p$ by a section of the root-finding algorithm, is indeed a minimum of $J$ along the trajectory; we verified that this is a general property for the family.}
    \label{fig:coords_along_traj}
\end{figure}

The time evolution of the Hill's coordinates along a single orbit of the family is shown in Fig. \ref{fig:coords_along_traj} for $J = 39^\circ$. We can check that the $J_p = J(t=0)$ defined in \eqref{eq:sigma1} is indeed the minimal mutual inclination (this has been checked numerically for all orbits). The initial condition on the orbit is fixed by \eqref{eq:sigma1}; choosing a different value for $w_1+w_2$ would yield a different result. For instance, if we had chosen $\sigma_1(\mathbf{x}) = w_1 + w_2 - \pi$ instead, we would have selected the maximal mutual inclination. However, as indicated above, along a trajectory, the variations of $J$ remain of order $\eps$ (see Fig. \ref{fig:coords_along_traj} bottom left). Therefore, choosing its maximum or minimum value as an initial condition has very little influence on the results presented below.  

\begin{figure}[htbp]
  \centering
  \includegraphics[width=0.9\textwidth]{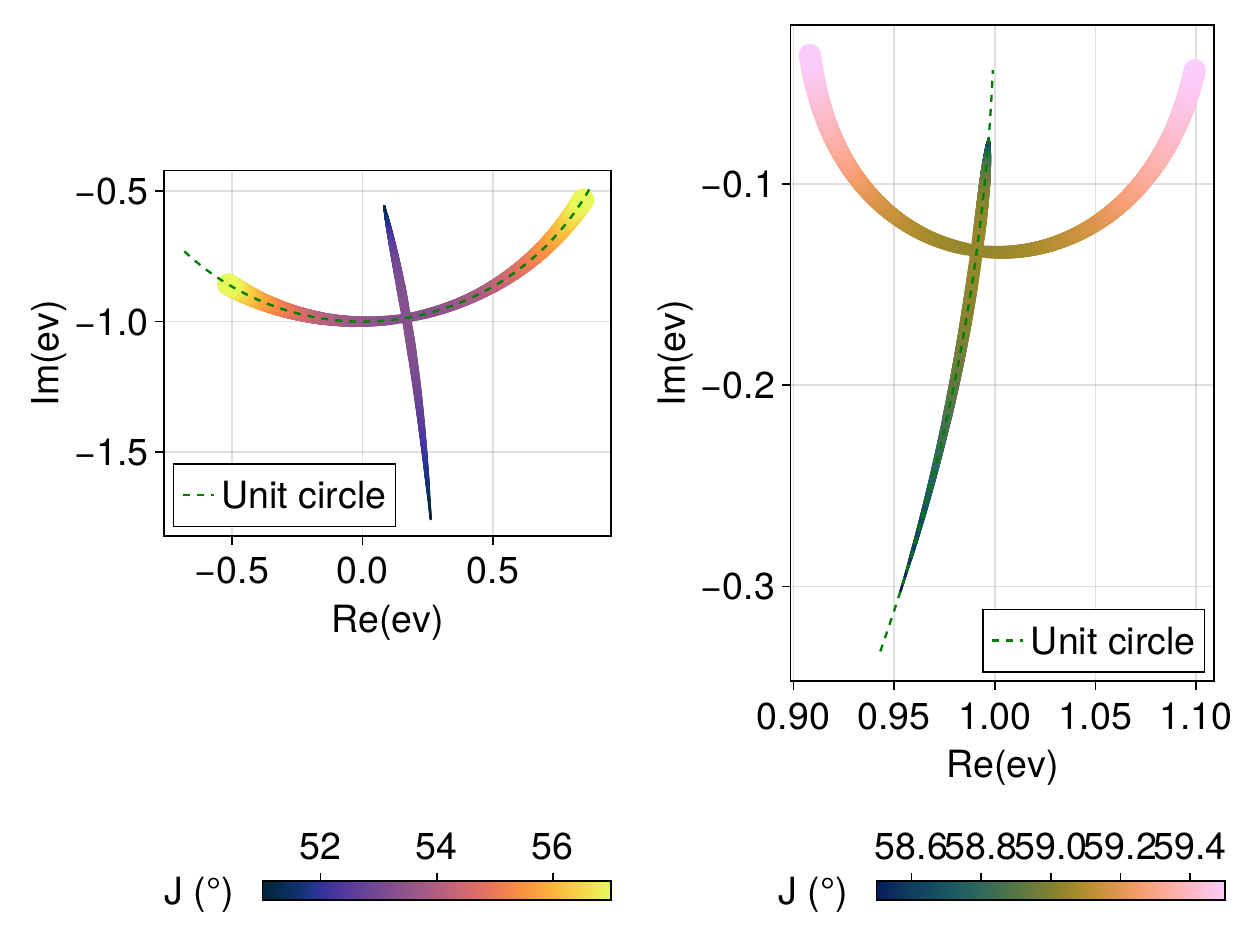}
  \caption{Trajectories of four monodromy matrix eigenvalues in the complex plane, along the stability changes of \VFL. Here, $\eps = 0.04$. The same collisions happen simultaneously for the complex conjugates, with positive imaginary parts. Curve width increases with $J$.}
  \label{fig:marchal_bifurc_eigenvalues}
\end{figure}

Similarly, we confirm that the planetary eccentricities (third row, right in Fig. \ref{fig:coords_along_traj}) are of order of the mass ratios, around $10^{-3}$. It is worth noting that all variables (apart from the non-periodic ones: $w_1 + w_2$ and $\Omega$) seem to actually be $T/2$-periodic. The longitude of the nodes $\Omega$ has a linear behavior associated with its precession,
 superposed with low-amplitude $T/2$-periodic oscillations. Note that, in the full problem, the precession frequency corresponding to this trajectory is  of $0.05542^{\circ}$ per period (period of $6496$ orbital revolutions), while we found $0.0555^{\circ}$ by application of the formula (\ref{eq:precession_noeud}). This indicates that, at least up to this value of the small parameter $\eps$, the average problem provides a good approximation of the trajectories of the full problem.  The analytical approximation (\ref{eq:prec_ana}) provides, for the same value of inclination, a frequency equal to $0.0523^{\circ}/\text{period}$, which suggests that this analytical approximation is still quite reasonable for $J=39^\circ$, even though it was obtained using an expansion at $J=0$.  
The family is stable from Lagrange ($J = 0^\circ$) up until about $60^\circ$; the branch starting from Euler, and everything above the $60^\circ$ threshold is unstable. Figure \ref{fig:marchal_bifurc_eigenvalues} shows the eigenvalue trajectories in the complex plane around the stable $\rightarrow$ unstable change, with four eigenvalues colliding in pairs and leaving the unit circle.\footnote{Two other bifurcations happen shortly after in the family: the four eigenvalues collapse on the real axis; then, two of them quickly fall back on the unit circle, colliding at $1$ and separating afterwards. We didn't show them, as they don't change the family's stability.} This stability change corroborates the frequency analysis results shown in figure \ref{fig:Moria} and \cite{Leleu2017PhD}; a long-term-stable region surrounds the linearly stable part of the family.

\subsubsection{\VFLtitre above Gascheau's masses}

Since the \VFL family appears to play a central role in the stability of the co-orbital regions surrounding equilateral equilibria, it is important to know whether its stability, at least partial (we have seen that this is the case for mutual inclinations ranging from $0^\circ$ to approximately $60^\circ$), persists for planetary masses higher than those considered so far, and to explore the consequences for the whole co-orbital regions.
Although Lagrange's equilateral (circular) configurations become unstable above the Gascheau value, here when $\eps \geq (3 - 2\sqrt2)/9 \approx 0.0191$, it has been shown that, in the elliptical case, stability can persist beyond this value \citep{Dan1964, Robe2002, Nauenberg2002}. It therefore seems relevant to ask whether a similar phenomenon exists for inclined quasi-circular orbits.

\begin{figure}[htbp]
  \centering
  \includegraphics[width=0.9\textwidth]{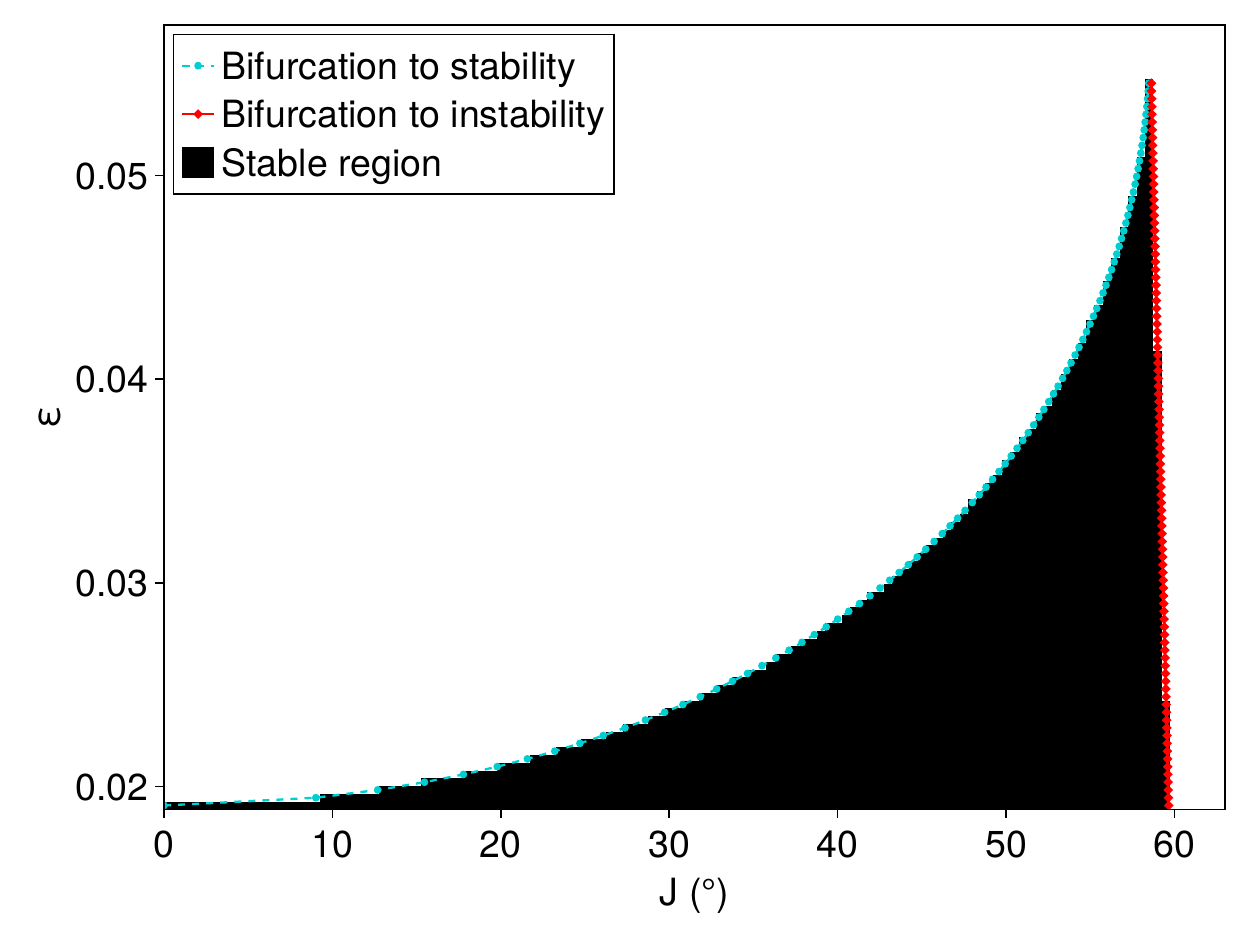}
  \caption{Stability of \VFL in the $(J, \eps)$ plane, with bifurcations up to $10^{-4}\phantom{1}^\circ$ on $J$.}
  \label{fig:marchal_chimney}
\end{figure}

In order to answer this question, we repeat the numerical analysis from section \ref{sec:num} for various mass parameters $\eps$; results are shown in figure \ref{fig:marchal_chimney}. Each pixel is a point of the family at a given mutual inclination $J_p$ and mass parameter $\eps$; its color is black if the orbit is linearly stable (the maximum distance of characteristic exponents to the unit circle is less than $10^{-6}$), and white otherwise. We confirm a stable region above Gascheau's value (which is at the figure's bottom): at a certain $J_p$ value (given by the green curve in Fig. \ref{fig:marchal_chimney}), the four eigenvalues from equation (\ref{eq:lagrange_eigenvalues}), which leave the unit circle above Gascheau, collide back in pairs, rejoining the circle along the family; the bifurcation is shown in the left plot of Fig. \ref{fig:marchal_bifurc_eigenvalues} for $\eps = 0.04$. The linear stability limit, where \VFL becomes unstable -- delimited by the red curve in fig \ref{fig:marchal_chimney} -- is the continuation of the bifurcation discussed in section \ref{vfl_planetary}; the mutual inclination $J_p$ where it happens slowly decreases with $\eps$. The stability remnant decreases in width with $\eps$, eventually disappearing for $\eps \approx 0.058$. 

\begin{figure}[htbp]
    \centering
    \includegraphics[width=0.9\linewidth]{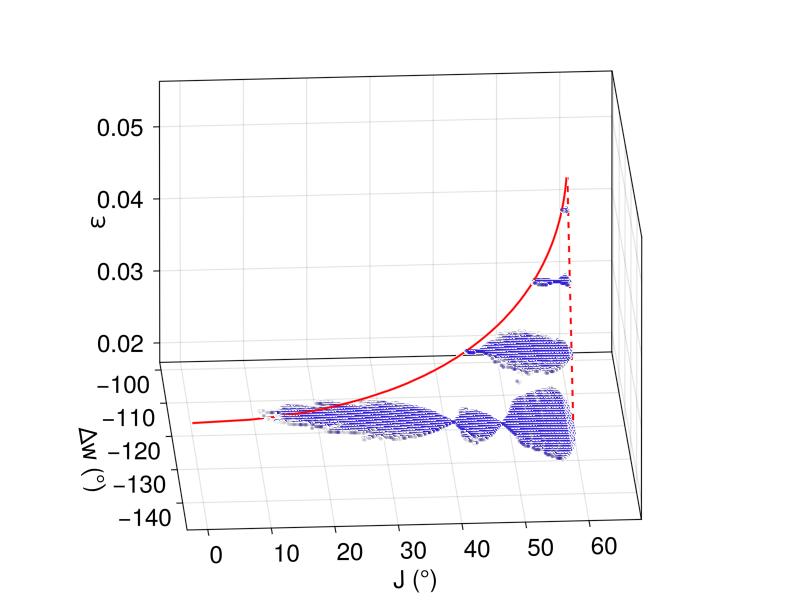}
    \caption{Stability maps (in blue) for $\eps = 0.02, 0.03, 0.04, 0.05$. The spectral stability limits are represented by red lines (see the text for more information).}
    \label{fig:chimney_3D}
\end{figure}

As we did in section \ref{sec:compar}, it is important to examine the possible links between the stability of \VFL and that of co-orbital regions for different values of planetary masses. For each value of $\eps$ in $\{0.02, 0.03, 0.04, 0.05\}$, we integrate, over $10\, 000$ orbital periods  in the full problem, a grid of initial conditions chosen as follows. We consider a discrete set $(r_j(J_p),R_j(J_p),w_j(J_p),G_j(J_p))$ of \VFL initial conditions parameterized by $J_p = k_J\delta_{J}$ with $k_J\in \{0,\cdots,320 \}$ and $\delta_J = 0.2^\circ$. In order to build a 2D initial conditions grid, we consider curves parallel to \VFL for which only $w_2$ is translated as $w_2 = w_2(J_2) + k_w\delta_w$ with $k_w\in \{-40,\cdots,40 \}$ and $\delta_w = 0.5^\circ$. The results of these simulations are plotted in \ref{fig:chimney_3D}, in the coordinates system $(J,\Delta w,\eps)$, where $\Delta w = w_1 - w_2$. The four blue horizontal layers correspond to regions of long-term stability that group together trajectories for which the relative diffusion in orbital frequency is less than $10^{-6}$, thus ensuring their stability \cite[see][]{RoLa2001}. As is the case for $\eps = 0.001$ (see Fig. \ref{fig:Moria}), the linearly stable section of the family \VFL occupies the center of these regions, whose boundary in the $J$ direction coincides with the red curves marking the stable-unstable transitions (solid lines for the lower bound and dotted lines for the upper bound) of \VFL reproduced in Fig. \ref{fig:marchal_bifurc_eigenvalues}.  As mentioned in section \ref{sec:compar}, the contours of stability domains are sculpted by secondary resonances, the most notable being the $2:1$ resonance between the libration frequency of the angle $w_1 - w_2$ and the mean motion. Note that the resonance $3:1$ between the same frequencies across the co-orbital zone for $\eps = 0.01$ (the figure is not reproduced here).

\section{Concluding remarks and beyond}

\cite{Marchal2009} introduced the vertical family we call \VFL in the average restricted three-body problem; we proved its existence beyond the restricted case, for small planetary masses of the average problem. We showcased numerical proof of its existence in the non-average problem, investigating its stability beyond planetary masses, for masses up to $1/8^\text{th}$ of the central masses.
Indeed, by analyzing the variational equation, we were able to identify certain critical values of planetary masses at which an increase in these masses at a fixed mutual inclination caused the associated periodic orbit normally unstable. For masses small enough to satisfy Gascheau's criterion, linear stability (planetary regime), even in the long term, remained guaranteed for all mutual inclinations less than a value close to $60^{\circ}$, a value that depends only slightly on the masses.  By increasing planetary masses, instability penetrates \VFL through weak inclinations until the stability zone disappears completely, the last region to disappear being close to $60^{\circ}$ for  $m_1  + m_2 \approx 0.058$, the sum of the three masses being equal to 1.
Although this critical value of mutual inclination had already been mentioned by certain authors (\cite{Leleu2017PhD} through numerical simulations,  \cite{QiRu2020} showing that tadpoles disappeared at around $61^\circ$ and \cite{GiLe2016} invoking a Lidov-Kozai mechanism that does not correspond to the situation explored in this work since the eccentricities are very low), its destabilization mechanism had not yet been clearly identified. 
Of course, the results described here are limited in the sense that we imposed the same mass on both small bodies. However, simulations not reported here show that the relevant parameter is not so much the individual masses as the sum of the planetary masses, a phenomenon already observed by \cite{LeRoCo2015}.  
In the exploration of the stability of \VFL mentioned by \cite{Marchal2009}, this study is limited  to low planetary masses. But continuing this work up to equal masses would bring definitive evidence of the link between \VFL and \pdouze{}. The latter being highly symmetrical, studying the symmetries of \VFL could also prove insightful.

\bmhead*{Data availability}
The code used for generating the figures in section \ref{sec:num} is publicly available as \cite{prieurSoftwareFindingInclined2025}.

\bmhead*{Acknowledgements}
We thank Jacques Fejoz for his precious insight and careful review.

\newpage
\begin{appendices}
\section{Expression of the average circular Hamiltonian using elliptic integrals}
\label{sec:appendix_Hamcirc}
%
%
Setting $p = \sin(\zeta_1/2)$ and $ u = p^2 +s_0^2 -p^2s_0^2$, the function $U_{s_0}$, defined in section (\ref{sec:circ_integrable}), can be rewritten as:
\be
\begin{aligned}
U_{s_0}(\zeta_1,\zeta_2) &= 2 - 2\cos\zeta_1 + 2s_0^2 \left(\cos\zeta_1 - \cos(\zeta_1 + 2\zeta_2) \right) \\
     & = 4\left[ p^2(1 - s_0^2) + s_0^2\sin^2(\zeta_1/2 + \zeta_2)\right]\\
     & = 4u\left[  1 - u^{-1}s_0^2\cos^2(\zeta_1/2 + \zeta_2)\right] \\
     & = 4u\left[  1 - u^{-1}s_0^2\sin^2(\zeta_1/2 + \zeta_2 - \pi/2)\right]. \\
\end{aligned}
\ee
It follows that:
\be
\begin{aligned}
\frac{1}{\pi}\int_0^\pi \frac{d\zeta_2}{\sqrt{U_{s_0}(\zeta_1,\zeta_2)}} 
&= \frac{1}{2\pi\sqrt u}\int_{\frac{\zeta_1-\pi}{2}}^{\frac{\zeta_1+\pi}{2}} \frac{d\theta}{\sqrt{1 - u^{-1}s_0^2\sin^2\theta}} \\
&=\frac{1}{2\pi\sqrt u} \left[ \cF\left(\frac{\zeta_1 +\pi}{2},\frac{s_0}{\sqrt u} \right) - \cF\left(\frac{\zeta_1 -\pi}{2},\frac{s_0}{\sqrt u} \right) \right],
\end{aligned}
\ee

where $\cF$ is the incomplete elliptic integral of first kind defined by:
\be
\cF(x,y) = \int_0^x \frac{dz}{\sqrt{1 - y^2\sin^2z}}.
\ee

 \section{Taylor expansion of the frequencies}
\label{sec:appendix_Taylor_freq}
In this appendix, the expansions in power series  of $s_0 = \sin(J_0/2)$ for the libration frequencies $\nu_L$ and the precession frequency of the node $\dot\Omega_p$ are given along \VFL. These expansions are limited to degree $20$, which provides a good approximation of the exact values for inclinations ranging from $0^0$ to approximately $60^\circ$.
According to (\ref{eq:fre_nu}), $\nu_L$ reads:
\be
\nu_L =  \sqrt{\eps}\ns\nut_L \sqrt{\frac{m_1+m_2}{m_0}} 
\ee
with 
\be
\begin{aligned}   
\nut_L\sqrt{\frac{4}{27}} =  &     1
    - \frac{4}{3}s_0^2
    + \frac{67}{36}s_0^4
    - \frac{23}{27}s_0^6
    - \frac{24721}{5184}s_0^8
    + \frac{14653}{3888}s_0^{10}
    + \frac{1165859}{62208}s_0^{12}\\
    &- \frac{914219}{46656}s_0^{14} 
    - \frac{3195892471}{35831808}s_0^{16}
    + \frac{8594100529}{80621568}s_0^{18}
    + \frac{1796301904997}{3869835264}s_0^{20}\\
    &+ \gO\left(s_0^{22}\right).
\end{aligned}
\ee

From (\ref{eq:precession_noeud}), the precession frequency of the ascending node is given by 
\be
\dot \Omega_p = \eps \ns\frac{m_1+m_2}{m_0}\sqrt{1 -4\gamma s_0^2} F_s, \qtext{with}\gamma = \frac{m_1m_2}{(m_1+m_2)^2}
\ee
and where $F_s$ satisfies the expansion: 
\be
\begin{aligned}   
F_{s} = &\phantom{-}\frac{3}{4}s_0^2
  - \frac{3}{2}\gamma s_0^4
 - \left( \frac{45}{32}
    + \frac{ 3}{2}\gamma^2\right)s_0^6
 + \left(   \frac{45}{16}\gamma
    - 3\gamma^3\right)s_0^8 \\
& + \left(   \frac{645}{128}
    + \frac{45}{16}\gamma^2
    - \frac{15}{2}\gamma^4\right)s_0^{10} 
 + \left( -\frac{ 645}{64}\gamma
    + \frac{45}{8}\gamma^3
    - 21\gamma^5\right)s_0^{12}\\
& + \left( - \frac{91325}{4096}
    - \frac{645}{64}\gamma^2
    + \frac{225}{16}\gamma^4
    - 63\gamma^6\right)s_0^{14} \\
 &+ \left(   \frac{91325}{2048}\gamma
    - \frac{645}{32}\gamma^3
    + \frac{315}{8}\gamma^5
    - \gamma^7\right)s_0^{16} \\
 &+ \left(   \frac{3602205}{32768}
    + \frac{91325}{2048}\gamma^2
    - \frac{3225}{64}\gamma^4
    + \frac{945}{8}\gamma^6
    - \frac{1287}{2}\gamma^8\right)s_0^{18}\\
 & + \left( - \frac{3602205}{16384}\gamma
    + \frac{91325}{1024}\gamma^3
    - \frac{4515}{32}\gamma^5
    + \frac{1485}{4}\gamma^7
    - 214\gamma^9\right)s_0^{20}\\
 &+ \gO\left(s_0^{22}\right).
\end{aligned} 
\label{eq:prec_ana}
\ee
\section{Lagrange equilateral relative equilibrium in Hill's coordinates}
\label{sec:appendix_Lagrange}

Although neither $w_j$ nor $\Omega_j$ are geometrically defined in the plane, Lagrange equilibria remain well defined in Hill coordinates.
 According (\ref{eq:L4_Hill}), all coordinates are constant except for true anomalies, which satisfy the differential equation
 \be
\dot w_j = \om\left(a + b\sin(2w_j+ (-1)^j\pi/3) + c\cos(2w_j+ (-1)^j\pi/3)\right),
\label{eq:edo_y}
\ee
where the coefficients $a,b$ and $c$ are given by:
\be
\begin{aligned}
 a &= (4m_0^4 + 6\sigma_0m_0^3 + (4\sigma_0^2 + p_0)m_0^2 + 5p_0\sigma_0m_0 + 2p_0^2)d^{-1}, \\
 b &= -\sqrt3p(m_1 - m_2)m_0 d^{-1}, \\
 c &= -p(3\sigma_0m_0 + 4m_0^2 + 2p_0) d^{-1},\\
 d &= m_0(2m_0 + m_1)(2m_0 + m_2)(m_0 + m_1 + m_2),\\
 p &= m_0m_1 + m_0m_2 + m_1m_2, \\ 
 p_0 &= m_1m_2, \quad \sigma_0 = m_1 + m_2, \\ 
 \end{aligned}
\label{eq:abc}
\ee
where a straightforward calculation leads the equality $a^2 - b^2 - c^2 = 1.$ 
The general solutions to this differential equation can be expressed as: 
$w_j(t) = w_j^0 + \omega t + F(\omega t) $, where $F$ is $2\pi$-periodic.

The dominant term of the 3 coefficients is given by:
\be
\begin{aligned}
a = & \, 1+ \eps^2a_2 +\gO(\eps^3)\qtext{with} a_2= \frac{(m_1 + m_2)^2}{2m_0^2} , \\
 b = & \,\eps^2b_2  +\gO(\eps^3) \qtext{with} b_2 = \frac{\sqrt3(m_2^2 - m_1^2)}{4m_0^2} ,\\
 c = & \, \eps c_1 + \eps^2c_2+\gO(\eps^3) \qtext{with} \\
 c_1 = &-\frac{(m_1 + m_2)}{m_0} \qtext{and} c_2= \frac{(3m_1^2 + 2m_1m_2 + 3m_2^2)}{4m_0^2}. \\
\end{aligned}
\ee
By neglecting terms in $\eps^3$ and more, the solution of (\ref{eq:edo_y}) reads:
\be
w_j(t) = w_j^0 + \om t + \eps\frac{c_1}{2}\sin\theta_j 
+\frac{\eps^2}{8}
\left( 
4c_2\sin\theta_j - 4b_2\cos\theta + c_1^2\sin2\theta_j
\right) + \gO(\eps^2)
\label{eq:w_Lagrange}
\ee
where $\theta_j = 2\omega t + (-1)^j\pi/3$.
%
%
As mentioned above, in the planar case, the angular position of the planets is identified by true longitude  $l_j = w_j + \Omega_j$. Since Lagrange's relative equilibrium is a periodic solution with frequency $\omega$,  the longitudes of the nodes $\Omega_j$ remain defined and are equal to: $\Om_j(t) = \Omega_j^0 - F(\omega t)$. These latter quantities therefore undergo only $\omega$-periodic oscillations.

\section{Difference between draconic mean anomalies in the average problem and draconic true anomalies in the full problem}
\label{sec:diff_avg_full}

\begin{figure}[!htb]
    \centering
    \includegraphics[width=0.85\linewidth]{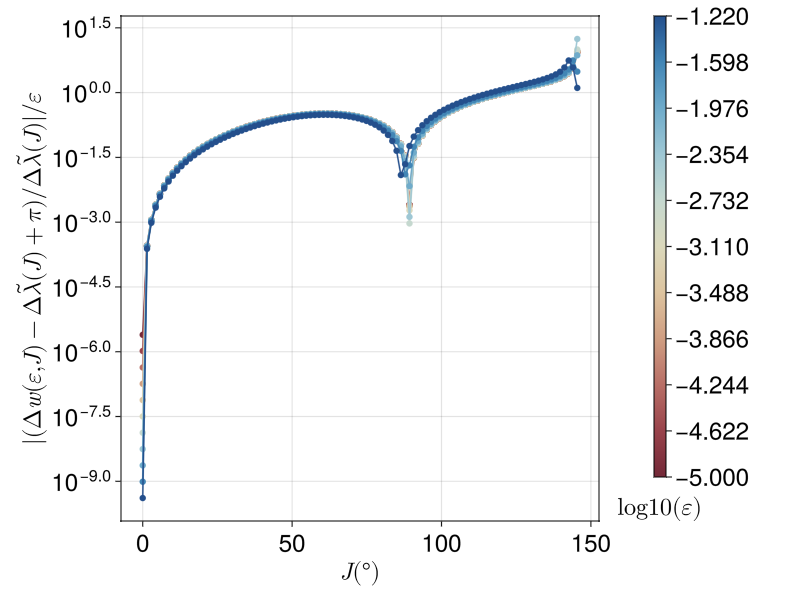}
    \caption{Relative variation of $\Delta \lam$ and $\Delta w$ along \VFL for various masses, divided by the mass parameter $\eps$. The addition of $\pi$ is due to the offset of the $\lat_j$. The drop near $J=100^\circ$ is due to the change of sign of the difference $\Delta w - \Delta \lam$.}
    \label{fig:comparison_avg_full}
\end{figure}

\changed{
Section \ref{sec:moyen} uses draconic mean anomalies, noted $\lat$ in the first-order circular average problem, whereas section \ref{sec:num} uses draconic true anomalies in the full problem. Figure \ref{fig:comparison_avg_full} shows the relative difference between $\Delta \lat$ and $\Delta w$, used as a family coordinate in figures (e.g. \ref{fig:Marchal} and \ref{fig:marchal_small_masses}), along \VFL, for mass parameters $\eps$ ranging from $10^{-5}$ to $0.06$. As Equation \eqref{eq:equil_circ} is independent of the masses, only $\Delta w$ is changing. The curves for the various masses are closely overlapping when divided by $\eps$; thus, the relative difference is mostly proportional to $\eps$, with a coefficient less than $1$, up to high mutual inclinations (around $J=120^\circ$). This gives an indication of the validity range of (\ref{eq:equil_circ}).
}

\section{\VFLtitre for unequal planetary masses}
\label{sec:appendix_unequal}
\begin{figure}[!htb]
     \centering
     \includegraphics[width=0.49\linewidth]{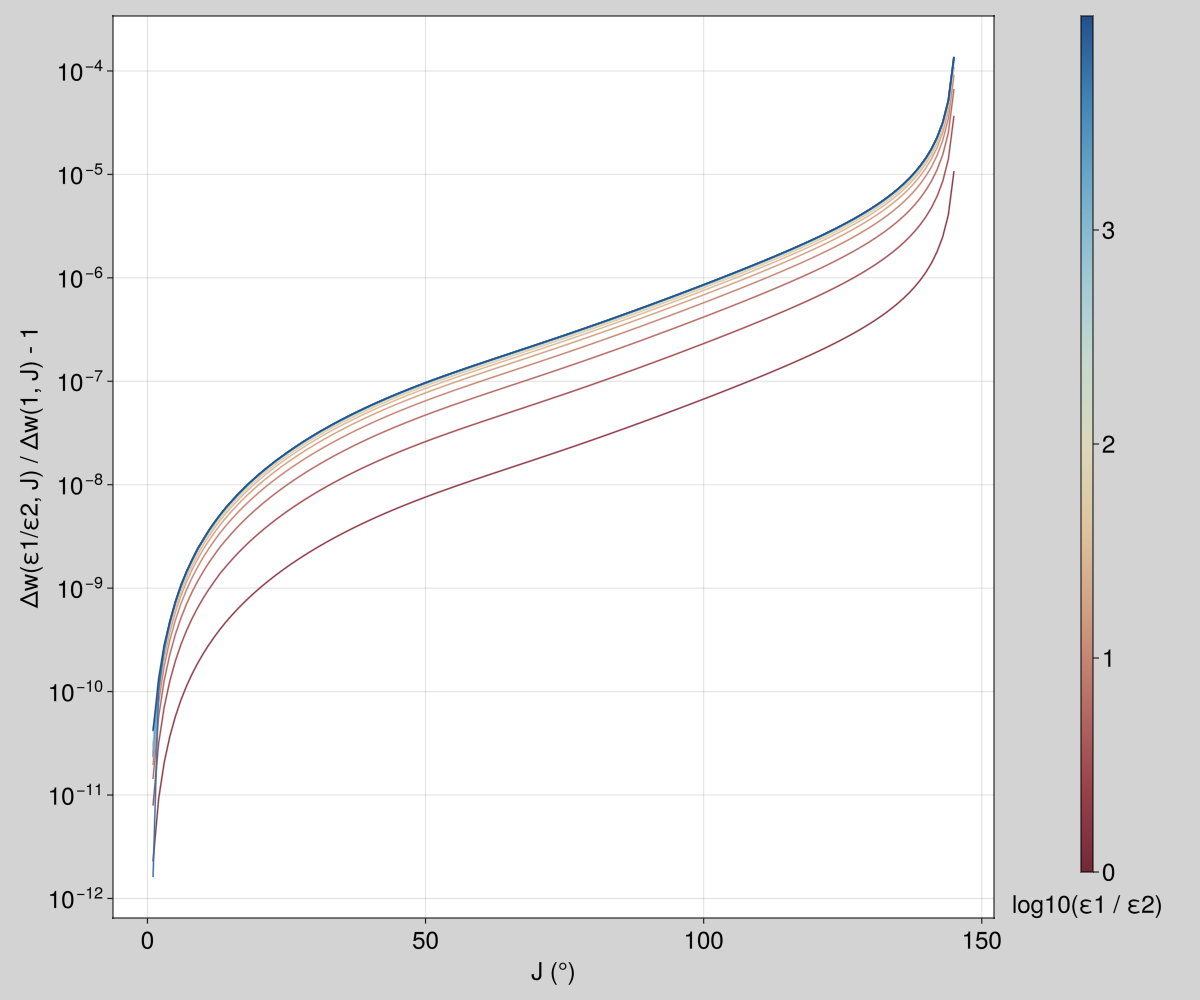}
     \includegraphics[width=0.49\linewidth]{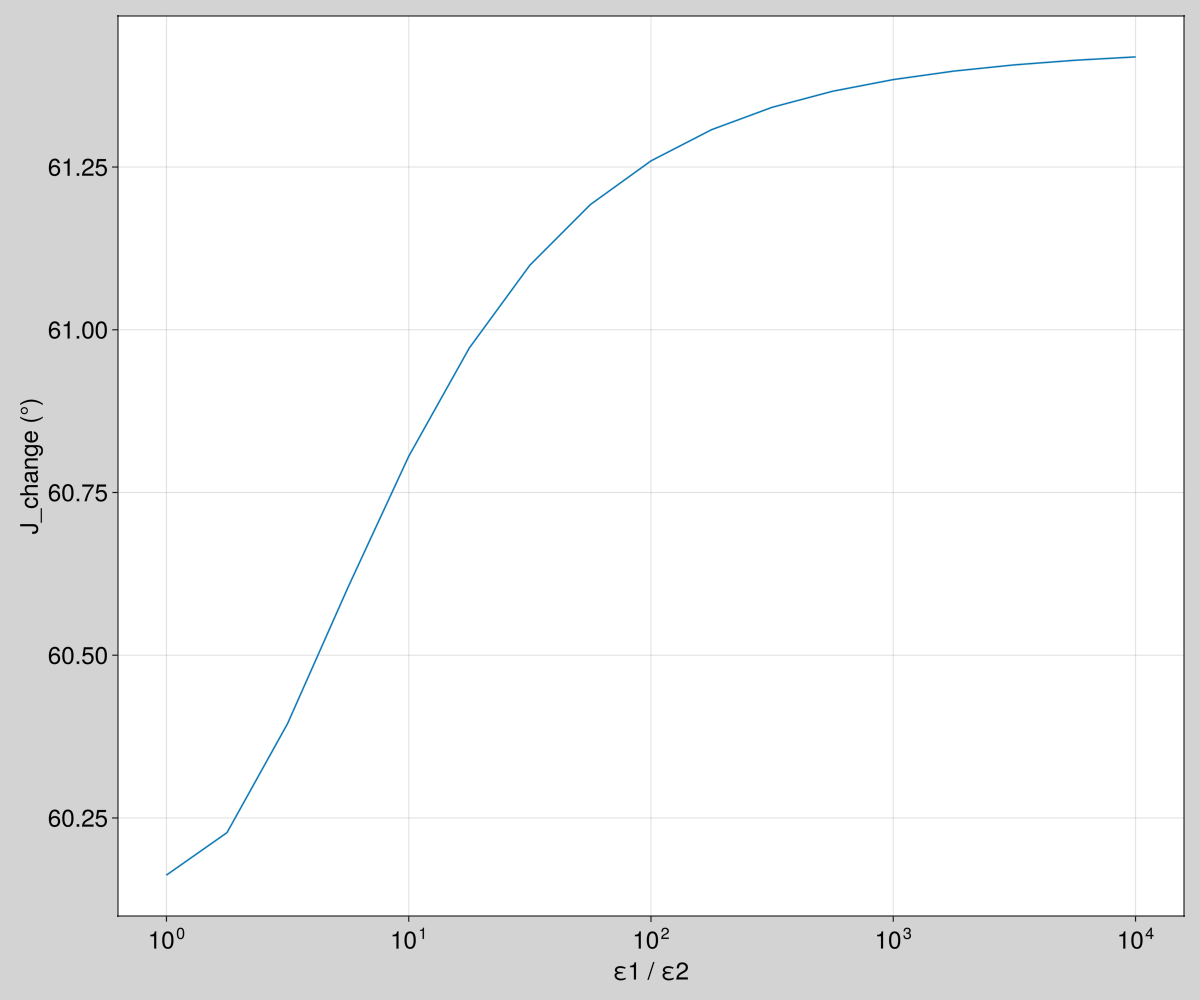}
     \caption{Left: Relative difference of $\Delta w$ between the unequal and equal masses case, for various aspect ratios, along \VFL. Right: Value of the mutual inclination for the stability change in \VFL for different mass ratios. In both cases, the sum of the two small masses is kept at $10^{-3}$, and the sum of all three masses at $1$.}
     \label{fig:unequal_masses}
 \end{figure}

\changed{We find \VFL for unequal planetary masses, with a ratio of small masses going up to $10^4$ -- keeping the sum of the two small masses at $10^{-3}$, and the sum of all three masses at $1$. Results are shown in Figure \ref{fig:unequal_masses} for the variation in $\Delta w$ and the position of the stability change, showing that qualitatively nothing is significantly changed, apart from the symmetry break between the two small masses. In both cases, results seem to converge to a limit, which one could imagine to be the restricted case, although as mentioned in Section \ref{sec:hill_coords}, we cannot attain it with the Jacobi reduction.}

\end{appendices}

\bibliography{papierHill1}

\end{document}